\documentclass{amsart}

\usepackage{amsmath}
\usepackage{amsthm}
\usepackage{amssymb}

\usepackage{cite}

\newtheorem{theorem}{Theorem}[section]
\newtheorem{lemma}[theorem]{Lemma}
\newtheorem{corollary}[theorem]{Corollary}
\newtheorem{proposition}[theorem]{Proposition}

\theoremstyle{definition}

\theoremstyle{remark}

\numberwithin{equation}{section}

\newcommand{\IR}{\mathbb{R}}
\newcommand{\IC}{\mathbb{C}}
\newcommand{\IZ}{\mathbb{Z}}
\newcommand{\IN}{\mathbb{N}}
\newcommand{\ID}{\mathbb{D}}
\newcommand{\IT}{\mathbb{T}}

\newcommand{\RE}{{\rm Re} \,}
\newcommand{\IM}{{\rm Im} \,}
\def\sech{\mathop{\rm sech}\nolimits}
\newcommand{\HS}{\mathcal{H}}
\newcommand{\BH}{\mathcal{B}(\mathcal{H})}

\begin{document}

\title{A Semigroup Composition C$^*$-algebra}

\author{Katie S. Quertermous}
\address{Katie S. Quertermous, Department of Mathematics, University of Virginia, Charlottesville, VA 22904}

\email{kgs5c@virginia.edu}
\subjclass[2000]{Primary  47B33; Secondary 47B20, 47B35, 47L80}
\thanks{\textit{Keywords:} Composition operator, Toeplitz operator, Hardy space, Almost periodic function, C$^*$-algebra, Commutator ideal}

\begin{abstract}
For $0 < s < 1,$ let $\varphi_s(z)=sz+(1-s).$  We investigate the unital C$^*$-algebra generated by the semigroup $\{C_{\varphi_s} : 0 < s < 1\}$ of composition operators acting on the Hardy space of the unit disk.  We determine the joint approximate point spectrum of a related collection of operators and show that the quotient of the C$^*$-algebra by its commutator ideal is isomorphic to the direct sum of $\IC$ and the algebra of almost periodic functions on the real line.  In addition, we show that the C$^*$-algebra is irreducible.
\end{abstract}

\date{September 4, 2009}

\maketitle

\section{\textsc{Introduction}}

For any analytic self-map $\varphi$ of the unit disk $\ID$, one can define the composition operator $C_{\varphi} : f \rightarrow f \circ \varphi,$ which is a bounded linear operator on the Hardy space $H^2(\ID).$  Individual composition operators on the Hardy space have been extensively studied, and there has been great success in relating the properties of a single composition operator $C_{\varphi}$ to the function-theoretic properties of the associated map $\varphi.$  Many of these results can be found in \cite{CowenMacCluer:1995} and \cite{Shapiro:1993}.

Recently, several authors have begun studying unital C$^*$-algebras generated by composition operators \cite{HamadaWatataniP, KrieteMacCluerMoorhouse:2007, KrieteMacCluerMoorhousePP, KrieteMacCluerMoorhousePP2, Jury:2007F, Jury:2007G, JuryUVA, JurySEAM}.  Although a few authors have considered composition operators induced by finite Blaschke products \cite{HamadaWatataniP, JuryUVA, JurySEAM}, most of the investigations have focused on composition operators induced by linear-fractional maps.  One motivation for this restriction is that the linear-fractional case has proven to be a useful model for guiding the study of more general composition operators in the single operator setting.  Moreover, composition operators induced by linear-fractional non-automorphisms serve as building blocks modulo the compact operators for certain more general composition operators \cite{KrieteMoorhouse:2007}.

The study of C$^*$-algebras generated by composition operators induced by linear-fractional maps tends to  split into two cases,  automorphism-induced generators and non-automorphism-induced generators.  The two cases have typically required different methods.
M. Jury has investigated the automorphism case \cite{Jury:2007F, Jury:2007G}, and Kriete, MacCluer and Moorhouse have studied the non-automorphism case \cite{KrieteMacCluerMoorhouse:2007, KrieteMacCluerMoorhousePP, KrieteMacCluerMoorhousePP2}.
The work of Kriete, MacCluer, and Moorhouse has focused on composition operators induced by maps $\varphi$ that either satisfy $\varphi(\zeta)=\eta$ for distinct points $\zeta$ and $\eta$ in the unit circle $\IT$   or fix a point $\zeta \in \IT$ and have $\varphi^{\prime}(\zeta)=1$.  In this paper, we begin consideration of the remaining non-automorphism case:  $\varphi(\zeta)=\zeta$ and $\varphi^{\prime}(\zeta) \neq 1.$

For $0 < s < 1,$ let $\varphi_s(z)=sz+(1-s).$  Note that $\varphi_s$ is a linear-fractional, non-automorphism self-map of $\ID$, $\varphi_s(1)=1$, and $\varphi^{\prime}_s(1)=s.$  For $0 < s, t <1,$ $C_{\varphi_s}C_{\varphi_t}=C_{\varphi_{st}}=C_{\varphi_t}C_{\varphi_s}$, so $\{C_{\varphi_s} : 0 < s < 1\}$ is a semigroup of commuting composition operators.  The elements of the semigroup have been  studied as individual operators by several authors, including Cowen and Ko \cite{CowenKoPP}, who determined the polar decomposition and Aluthge transform of $C_{\varphi_s}$, and Cowen and Kriete \cite{CowenKriete:1988}, who showed that $C_{\varphi_s}^*$ is subnormal.  In this paper, we determine the structure of the unital C$^*$-algebra generated by the semigroup modulo its commutator ideal.  Recall that the commutator ideal of a C$^*$-algebra $\mathcal{A}$ is the closed ideal of $\mathcal{A}$ generated by elements of the form $[A,B]=AB-BA$ for $A, B \in \mathcal{A}.$

To simplify our investigations, we will consider several sets of operators that are unitarily equivalent to $\{C_{\varphi_s} : 0 < s <1\}.$  These alternate settings are more conducive to the determination of spectral information and invariant subspaces than the original setting in $H^2(\ID).$   The operators that we will consider, and the Hilbert spaces on which they act, are introduced in Section 2.

In Section 3, we recall the definition of the joint approximate point spectrum of a collection of operators and a theorem of J. Bunce that relates the structure of the unital C$^*$-algebra generated by a collection of commuting hyponormal operators to the joint approximate point spectrum of the collection.  Motivated by this relationship, we determine the joint approximate point spectra of finite subsets of $\{T_{s^z} : 0 < s< 1\}$, a collection of operators that is unitarily equivalent to $\{C_{\varphi_s}^* : 0 < s < 1\}.$  These spectra prove to be rather complicated, which makes a direct application of Bunce's theorem impractical.  Instead we connect these spectra to the joint approximate point spectra of almost periodic Toeplitz operators, a better-understood class of operators.  In Section 4, we use this connection to prove our main result:
\begin{theorem} 
For $0 < s <1,$ let $\varphi_{s}(z)=sz+(1-s).$  Let $\mathcal{C}$ denote the commutator ideal of $C^*\left(\left\{C_{\varphi_s}: 0 < s < 1\right\}\right).$  Then there exists a $*$-homomorphism $\psi :C^*\left(\left\{C_{\varphi_s}: 0 < s < 1\right\}\right) \rightarrow  AP(\IR) \oplus \IC$ such that    $$ 0 \rightarrow \mathcal{C}  \hookrightarrow C^*\left(\left\{C_{\varphi_s} :
   0 < s < 1\right\} \right) \stackrel{\psi}{\rightarrow} AP(\mathbb{R}) \oplus \mathbb{C} \rightarrow 0$$ is a short exact sequence.
 \end{theorem}
 \noindent Here, and throughout this paper, $C^*\left(\left\{C_{\varphi_s} : 0 < s < 1\right\}\right)$ denotes the unital C$^*$-algebra generated by $\{C_{\varphi_s} : 0 < s < 1\}.$

  In Section 5, we show that the commutant of  $C^*\left(\left\{C_{\varphi_s} : 0 < s < 1\right\}\right)$ consists of only scalar multiples of the identity, so  $C^*\left(\left\{C_{\varphi_s} : 0 < s < 1\right\}\right)$ is irreducible.

\section{\textsc{Some Related Spaces and Operators}}

During the course of our investigations, we will use a variety of spaces and operators.  We now describe these spaces and set up our notation.
\subsection{Hardy space of the Disk}

The Hardy space of the disk, denoted $H^2(\ID),$ is the set of all functions $f(z)=\sum_{n=0}^{\infty}a_nz^n$ analytic in the open unit disk $\ID$ that satisfy $$||f||_{H^2(\ID)}^{2} := \sum_{n=0}^{\infty} |a_n|^2 < \infty.$$  The Hardy space has reproducing kernels $\kappa_w = (1-\overline{w}z)^{-1}$ for all $w \in \ID$ that satisfy $\left<f, \kappa_w\right>_{H^2(\ID)}=f(w)$ for all $f \in H^2(\ID).$  More information about $H^2(\ID)$ can be found in \cite{CowenMacCluer:1995} and \cite{Duren:1970}.

\subsection{A Unitarily Equivalent Space: $H^2(\mu)$}
We define a measure $\mu$ on the half-plane $\Omega=\left\{z \in \IC : \RE z \geq -\frac{1}{2}\right\}$ by \begin{equation}\label{measure}d\mu = \sum_{n=-1}^{\infty} 
\frac{\left|\Gamma\left(\frac{n}{2}+iy+1\right)\right|^2}{2\pi(n+1)!}dy d\delta_{\frac{n}{2}}(x),\end{equation} where $\delta_{\frac{n}{2}}$ is the measure on $\IR$ having a unit point mass at $x=\frac{n}{2}$ and $\Gamma$ is the gamma function.  The measure $\mu$ is finite with total mass equal to 1 \cite{KrieteTrutt:1971}. 
For convenience, we denote $L^2\left(\Omega, \mu\right)$ and $L^{\infty}\left(\Omega, \mu\right)$  by $L^2(\mu)$ and $L^{\infty}(\mu),$ respectively.  For $0< a \leq 1,$ the function $f(z)=a^z$ is in  $L^2(\mu)$. The Hilbert space $H^2(\mu)$ is defined as the closed linear span of $\{a^z: 0<a \leq 1\}$ in $L^2(\mu)$.  

The orthogonal projection of $L^2(\mu)$ onto $H^2(\mu)$ will be denoted $P_{\mu}.$  For $f \in L^{\infty}(\mu),$ we define the multiplication operator $M_{f}: L^2(\mu) \rightarrow L^2(\mu)$ by $M_{f}\,g=fg$ for all $g \in L^2(\mu)$ and the Toeplitz-like operator $T_{f}: H^{2}(\mu) \rightarrow H^2(\mu)$ by $T_f\,h=P_{\mu}fh$ for all $h \in H^2(\mu).$  Note that for $0< s, t <1,$ $T_{s^z}T_{t^z} = T_{(st)^z}=T_{t^z}T_{s^z}$, and $T_{s^z}=M_{s^z} |_{H^2(\mu)}$ since $s^za^z \in H^2(\mu)$ for all $a \in (0,1]$. Thus, $\{T_{s^z}: 0<s <1\}$ is a semigroup of commuting subnormal operators.  

A unitary operator $V : H^2(\ID) \rightarrow H^2(\mu)$ can be defined in the following way \cite{CowenKriete:1988}: 
 For $0< a, b \leq 1$, \begin{equation}\label{preserve} \left<\kappa_{(1-a)}, \kappa_{(1-b)}\right>_{H^2(\ID)} = \left<a^z, b^z\right>_{H^2(\mu)}. \end{equation}  
 Define the operator $V$ on the functions $\kappa_{(1-a)}$ by $V\kappa_{(1-a)} = a^z$.  The set $\{\kappa_{1-a}: 0 < a \leq 1\}$ has dense linear span in $H^2(\ID)$, and the set $\{a^z : 0 < a \leq 1\}$ has dense linear span in $H^2(\mu).$  Since $V$ preserves inner products by (\ref{preserve}), it has a unique extension to a unitary operator, also denoted $V$, from $H^2(\ID)$ onto $H^2(\mu).$  
 
 Our interest in the space $H^2(\mu)$ comes from the following theorem of Cowen and Kriete:

\begin{theorem}\label{cowenkrieteunitary}{\rm \cite[Theorem 18]{CowenKriete:1988}} Let $0 < s <1.$  Then $C_{\varphi_s}^*$ is unitarily equivalent to $T_{s^z}$ via $V.$ \end{theorem}

We will find that for many of our investigations it is more convenient to study $C^*(\{T_{s^z} : 0 < s < 1\})$ instead of considering $C^*(\{C_{\varphi_s}: 0 < s <1\})$ directly.

\subsection{A Second Equivalent Space: The Newton Space}
 
 The Newton space, denoted $\mathcal{N}$, is a Hilbert space of analytic functions on $\left\{z \in \IC : \\ \RE z > - \frac{1}{2}\right\}$ that  has the Newton polynomials,  
 \begin{equation*}  N_n(z) = \left\{\begin{array}{cc} 1, & n =0 \\ (-1)^n \frac{z(z-1) \ldots (z-(n-1))}{n!}, & n \geq 1\end{array} \right.,\end{equation*} as an orthonormal basis.  Recall that the Newton polynomials satisfy \begin{equation}\label{Newtonsum} (1-w)^z = \sum_{n=0}^{\infty} N_n(z) w^n,\end{equation} for $|w| < 1$ and $z \in \IC$.
The reproducing kernels for $\mathcal{N}$ are \begin{equation}\label{NewtonRK} K_w(z)=\frac{\Gamma(z+\overline{w}+1)}{\Gamma(z+1)\Gamma(\overline{w}+1)} \end{equation} for $ w \in \IC$ with $\RE w > -\frac{1}{2}.$

   Every function $f \in \mathcal{N}$ has a non-tangential boundary function $f(\zeta)$ that is defined a.e. on the line $\RE \zeta = -\frac{1}{2}$, and $\mathcal{N}$ is contained isometrically in $L^2(\mu)$ by taking the values of a function on the line $\RE z = -\frac{1}{2}$ to be those of the boundary function \cite{MarkettRosenblumRovnyak:1986}.  In fact, since the linear span of 
  $\{a^z : 0 < a \leq 1\}$ is dense in $\mathcal{N}$ \cite{KrieteTrutt:1971}, we can view $\mathcal{N}$ as a subspace of $H^2(\mu).$ 
Moreover, the map $U_1$ that restricts a function $g \in H^2(\mu)$ to the representative function $f(w)=\left<g, K_w\right>_{H^2(\mu)}$ on the open half-plane $\left\{w : \RE w > -\frac{1}{2}\right\}$ is a unitary operator from $H^2(\mu)$ onto $\mathcal{N}$ that satisfies $f=g$ $\mu$-almost everywhere on the open half-plane.  Note that $U_1 a^z=a^z$ for all $0 < a \leq 1.$ 

\subsection{A Third Equivalent Space: $H^2(\tilde{\mu})$}

The map $\Psi(z)=\frac{z}{1-z}$ takes $\overline{\ID}$ onto $\Omega.$  We define the measure $\tilde{\mu}$ of $\overline{\ID}$ by $\tilde{\mu}=\mu \circ \Psi$ and $\tilde{\mu}(\{1\})=0.$  Let $L^2(\tilde{\mu})$ and $L^{\infty}(\tilde{\mu})$ denote $L^{2}(\overline{\ID}, \tilde{\mu})$ and $L^{\infty}(\overline{\ID}, \tilde{\mu})$, respectively.  We define the Hilbert space $H^2(\tilde{\mu})$ to be the closed linear span of  $\left\{a^{\frac{z}{1-z}} : 0 < a \leq 1\right\}$ in $L^2(\tilde{\mu}).$  It is clear that the map $U_2$, defined by $$U_2\left(\sum_{j=1}^n a_j^z\right)=\sum_{j=1}^n a_j^{\frac{z}{1-z}}$$ for $n \in \IN$ and $0 < a_j \leq 1$ for $j=1, \ldots, n$, extends to a unitary operator from $H^2(\mu)$ onto $H^2(\tilde{\mu}).$

We denote the projection of $L^2(\tilde{\mu})$ onto $H^2(\tilde{\mu})$  by $P_{\tilde{\mu}}$.
For $f \in L^{\infty}(\tilde{\mu}),$  we define the multiplication operator $M_{f}: L^2(\tilde{\mu}) \rightarrow L^2(\tilde{\mu})$ by $M_{f}\,g=fg$ for all $g \in L^2(\tilde{\mu})$ and the
 Toeplitz-like operator $\tilde{T}_{f}: H^{2}(\tilde{\mu}) \rightarrow H^2(\tilde{\mu})$ by $\tilde{T}_f\,h=P_{\tilde{\mu}}fh$ for all $h \in H^2(\tilde{\mu}).$  Note that for $0< s<1,$ $\tilde{T}_{s^{\frac{z}{1-z}}}=M_{s^{\frac{z}{1-z}}} |_{H^2(\tilde{\mu})}$ and $U_2T_{s^z} = \tilde{T}_{s^{\frac{z}{1-z}}}U_2.$
 
One advantage of studying operators in this space is that the support of $\tilde{\mu}$ is the union of the circles $C_m$ with centers $\frac{m}{m+1}$ and radii $\frac{1}{m+1}$ for $m=0,1, 2, \ldots,$ a compact set in $\IC.$  Also, $H^2(\tilde{\mu})$ is the closure of the polynomials in $L^2(\tilde{\mu})$  \cite{KrieteTrutt:1971}.  These facts will be key ingredients in  proving that $C^{*}\left(\left\{C_{\varphi_s} : 0 < s < 1\right\}\right)$ is irreducible.

\subsection{Hardy Space of a Strip}
For a set $J \subset \IC$, let $Hol\left(J\right)$ denote the set of all functions analytic on $J$.
For $a, b \in \IR$ with $a <b$,  
we define a strip $S(a,b) :=\{z \in \IC : a < \RE z <b$\}.  
Then the Hardy class $H^2\left(a,b\right)$ for the strip $S(a,b)$ is \begin{equation*} \left\{ F \in Hol\left(S(a,b)\right) : \sup_{a < x< b}\int_{\IR} |F\left(x+iy\right)|^2 dy < \infty\right\}. \end{equation*}  Functions $F(z)$ in $H^2(a,b)$ have boundary functions  \begin{equation*} F(a+iy) = \lim_{x \rightarrow a^+} F(x+iy), \, \, \, \, F(b+iy)=\lim_{x \rightarrow b^{-}} F(x+iy), \end{equation*} 
 which exist almost everywhere and in the metric of $L^2(\IR)$\cite{RosenblumRovnyak:1994} and satisfy a Three Lines-type Lemma:

\begin{theorem}\label{three lines} {\rm \cite[Theorem 2.3]{BakanKaijser:2007}} Let $F \in H^2(a,b).$  Let $a \leq \alpha < \beta < \gamma \leq b.$ Then \begin{equation*} \left|\left|F\left(\beta + i y\right)\right|\right|_{L^2\left(\IR\right)} \leq \left|\left|F\left(\alpha + i y \right)\right|\right|^{\frac{\gamma-\beta}{\gamma-\alpha}}_{L^2\left(\IR\right)} \left|\left|F\left(\gamma +  i y \right)\right|\right|^{\frac{\beta-\alpha}{\gamma-\alpha}}_{L^2\left(\IR\right)}. \end{equation*} \end{theorem}

The Newton space is related to the Hardy space of a strip by the following theorem:

\begin{theorem}\label{MRRThm} {\rm \cite[Theorem 1.6]{MarkettRosenblumRovnyak:1986}} If $f\left(z\right)\in \mathcal{N}$, then $\Gamma\left(z+1\right)f\left(z\right)$ is in  $H^2\left(-\frac{1}{2}, \xi\right)$ for every $\xi \in \left(-\frac{1}{2}, \infty\right).$ \end{theorem}

\subsection{Hardy Space of the Line and Almost Periodic Functions}
We set $L^{2}(\IR):=L^2(\IR, m)$, where $m$ is
Lebesgue measure, and we denote the Fourier transform on $L^2(\IR)$ by $\mathcal{F}$.  The Hardy space $H^2(\IR)$ is the subspace of $L^2(\IR)$ of all functions $f$ for which $\mathcal{F}f$ is supported on $\IR^+$.  The subspace is non-trivial and is the subspace of $L^2(\IR)$ of boundary values of functions in the Hardy space of the upper half-plane.

We are also interested in a class of continuous functions on $\IR$. For $\alpha \in \IR$, we define $\chi_{\alpha}: \IR \rightarrow \IC$ by $\chi_{\alpha}(y)=e^{i\alpha y}$ for all $y\in \IR$.  Finite linear combinations of the functions $\{\chi_{\alpha} : \alpha \in \IR\}$ are called trigonometric polynomials.  A continuous function $f: \IR \rightarrow \IC$ is  called almost periodic if, for all $\varepsilon >0$, there exists a trigonometric polynomial $T_{\varepsilon}(x)$ such that $|f(x)-T_{\varepsilon}(x)| < \varepsilon$ for all $x \in \IR.$  The space of all almost periodic functions is denoted $AP(\IR).$  The theory of almost periodic functions was created by H. Bohr and has been thoroughly developed over the course of the last century \cite{Bohr:1947, Corduneanu:1968, Besicovitch:1955, LevitanZhikov:1982}.

The orthogonal projection of 
$L^2(\IR)$ onto $H^2(\IR)$ will be denoted by $P_m$.
For $f \in L^{\infty}(\IR):=L^{\infty}(\IR, m)$, we define  the multiplication operator $M_{f}: L^2(\IR) \rightarrow L^2(\IR)$ by $M_{f}g =  fg$ for all $g \in L^2(\IR)$ and the Toeplitz operator $W_{f}: H^2(\IR) \rightarrow H^2(\IR)$ by $W_{f}h=P_m fh$ for all $h \in H^2(\IR)$.   We are particularly interested in the collection $\left\{W_{\chi_{\alpha}} : \alpha \in \IR^{+}\right\}$.  If $\alpha \in \IR^{+}$, then $W_{\chi_{\alpha}}=M_{\chi_{\alpha}}|_{H^2(\IR)}.$  The C$^*$-algebra $C^*\left(\left\{W_{\chi_{\alpha}}: \alpha  \in \IR^+\right\}\right)$ has been studied extensively \cite{CoburnDouglas:1969, CoburnDouglas:1971, CoburnDouglasSchaefferSinger:1971}, and its structure modulo its commutator ideal is described by the following result:

\begin{theorem}\label{CDiso}{\rm \cite{CoburnDouglas:1971}} Let $\mathcal{C_W}$ be the commutator ideal of $C^*\left(\left\{W_{\chi_{\alpha}}: \alpha  \in \IR^+\right\}\right)$.  Then \begin{equation*} C^*\left(\left\{W_{\chi_{\alpha}}: \alpha  \in \IR^+\right\}\right) / \mathcal{C_W} \cong AP(\IR).\end{equation*} \end{theorem}

\section{\textsc{Joint Approximate Point Spectra}}


The joint approximate point spectrum, denoted $\sigma_{ap}(A_1, \ldots, A_n)$, of a finite set $\{A_1, \ldots, A_n\}$ of commuting bounded 
operators on a Hilbert space $\mathcal{H}$ is the set \begin{equation*}\{(\lambda_1, \ldots, \lambda_n) \in \IC^n:
\BH(A_1-\lambda_1 I) + \ldots + \BH(A_n-\lambda_nI) \neq \BH\}.\end{equation*} By
the work of J. Bunce \cite{Bunce:1971}, $\sigma_{ap}(A_1, \ldots, A_n)$ is a
non-empty compact set, and \begin{equation}\label{inclusionindirectsum}\sigma_{ap}(A_1,  \ldots, A_n)
\subset \sigma_{ap}(A_1) \times \sigma_{ap}(A_2) \times \ldots \times
\sigma_{ap}(A_n).\end{equation}  An equivalent
characterization is that
$(\lambda_1, \ldots, \lambda_n) \in \sigma_{ap}(A_1,  \ldots,
A_n)$ if and only if there exists a sequence $\{x_m\}$ of unit
vectors in $\HS$ such that \begin{equation*} \lim_{m \rightarrow \infty} ||(A_j-\lambda_j I)x_m||_{\HS} =0\end{equation*} for all $j\in \{1, \ldots, n\}$ \cite{EFT:1972}.  Thus, for $n =1$, $\sigma_{ap}(A_1)$ is the usual approximate point spectrum, i.e.\ the set of all $\lambda \in \IC$ such that $A_1-\lambda I $ is not bounded below.

The connection between the joint approximate point spectrum of a collection of hyponormal operators and the C$^*$-algebra generated by the operators is identified in the following theorem of Bunce:

 \begin{theorem}\label{Buncefinite}{\rm \cite[Corollary 4]{Bunce:1971}}  If $A_1, A_2, \ldots, A_n$ are commuting hyponormal operators, then $ \sigma_{ap}(A_1,  \ldots, A_n)$ equals \begin{equation*} \left\{\left(\rho(A_1),  \ldots, \rho(A_n)\right) : \rho \, \text{is a character on} \, \, C^*(\{A_1,  \ldots, A_n\})\right\}, \end{equation*} and if \begin{equation*} J= \bigcap \left\{\rho^{-1}(0) : \rho \, \text{is a character on} \, \, C^*(\{A_1,  \ldots, A_n\})\right\}, \end{equation*} then \begin{equation*} C^{*}\left(\left\{A_1, \ldots, A_n\right\}\right) / J \cong C(\sigma_{ap}(A_1, \ldots, A_n)). \end{equation*} \end{theorem}
 
 Note that the ideal $J$ is the commutator ideal of $C^*(\{A_1, \ldots, A_n\}),$ and the map from $C^{*}\left(\left\{A_1, \ldots, A_n\right\}\right) / J$ onto $C(\sigma_{ap}(A_1, \ldots, A_n))$ is simply the Gelfand transform.  
 
 Both the definition of the joint approximate point spectrum and Theorem \ref{Buncefinite} can be extended to infinite collections of operators. For a family $\mathcal{S}=\{A_{\alpha} : \alpha \in \Lambda\}$ of commuting hyponormal operators, we define \begin{equation*} \sigma_{ap}(\mathcal{S})=\left\{\{\rho(A_{\alpha})\}_{\alpha \in \Lambda} : \rho \, \, \text{is a character on} \, \, C^*(\mathcal{S})\right\}. \end{equation*}  By Proposition 5 in \cite{Bunce:1971}, $\sigma_{ap}(\mathcal{S})$ is the inverse limit of the sets $\sigma_{ap}(A_{\alpha}: \alpha \in F),$ where  $F \subset \Lambda$ is finite, and thus $\sigma_{ap}(\mathcal{S})$ is a compact set.  As suggested by the notation, we call $\sigma_{ap}(\mathcal{S})$ the joint approximate point spectrum of $\mathcal{S}$. 
 
\begin{theorem}\label{Bunceinfinite}{\rm \cite{Bunce:1971}}  Let $\mathcal{S}=\{A_{\alpha} : \alpha \in \Lambda\}$ be a family of commuting hyponormal operators.  Let $\mathcal{C}$ be the commutator ideal of $C^*(\mathcal{S}).$  Then \begin{equation*} C^*(\mathcal{S}) / \mathcal{C} \cong C(\sigma_{ap}(\mathcal{S})). \end{equation*} \end{theorem}

\subsection{Calculating $\sigma_{ap}(T_{s_1^z}, \ldots, T_{s_n^z})$}

Since $\{T_{s^z} : 0 < s < 1\}$ is a collection of commuting subnormal, and hence hyponormal, operators,  Theorems \ref{Buncefinite} and \ref{Bunceinfinite} can be applied to the unital C$^*$-algebras generated by these operators.  Thus, we wish to determine the joint approximate point spectra of finite subsets of $\{T_{s^z}  : 0 < s < 1\}$.  We begin our investigations with an inner product calculation.

\begin{lemma}\label{innerprodcalc} Let $0< s,t <1$ and $w \in \IC$ with $\RE w > -\frac{1}{2}.$  Let $k_w=\frac{K_w}{||K_w||_{\mathcal{N}}}$ be the normalized reproducing kernel function for $\mathcal{N}$ corresponding to evaluation at $w$.  Then \begin{equation*} \left<T_{t^z}^* T_{s^z}k_w, k_w\right>_{H^2(\mu)}
= \frac{s^wt^{\overline{w}}}{(s+t-st)^{2 \RE w +1}}. \end{equation*} \end{lemma}

\begin{proof}

Recall from (\ref{NewtonRK}) that \begin{equation*} K_w(z)=\frac{\Gamma(z+\overline{w}+1)}{\Gamma(z+1)\Gamma(\overline{w}+1)} \, \, \, \text{and} \, \, \, 
 ||K_w||^2_{H^2(\mu)} = ||K_w||^2_{\mathcal{N}} = \frac{\Gamma(2 \RE w +1)}{|\Gamma(\overline{w}+1)|^2}. \end{equation*}  Thus, we find that 
\begin{align}   &\left<T_{t^z}^*T_{s^z}k_w,  k_w\right>_{H^2(\mu)} \notag\\ & \qquad \qquad=  \left<M_{s^z}k_w, M_{t^z} k_w \right>_{L^2(\mu)}  \nonumber \\
& \qquad \qquad = \int_{\Omega} s^z t^{\overline{z}} \frac{|\Gamma(z+\overline{w}+1)|^2}{|\Gamma(z+1)|^2|\Gamma(\overline{w}+1)|^2}\frac{|\Gamma(\overline{w}+1)|^2}{\Gamma(2 \RE w+1)} d \mu \nonumber  \\
& \qquad \qquad = \sum_{n=-1}^{\infty} \int_{-\infty}^{\infty} s^{\frac{n}{2}+iy}t^{\frac{n}{2}-iy}\frac{\left|\Gamma\left(\frac{n}{2}+iy+\overline{w}+1\right)\right|^2 dy}{\Gamma(2\RE w +1)2\pi(n+1)!}  \nonumber \\
& \qquad \qquad = \sum_{n=-1}^{\infty} \frac{(st)^{\frac{n}{2}}\left(\frac{s}{t}\right)^{i \IM w}}{2\pi(n+1)!}\int_{-\infty}^{\infty} \frac{e^{-i\ln\left(\frac{t}{s}\right)\alpha}\left|\Gamma\left(\frac{n+2+2 \RE w}{2}+i\alpha\right)\right|^2 d \alpha}{\Gamma(2\RE w +1)} \label{ipsofar}
\end{align}
by applying the definition of $\mu$ and using the change of variable $\alpha=y-\IM w$ to obtain the last line.  Note that, 
for $n \in \{-1, 0, 1, \ldots\}$, \begin{equation}\label{SwitchingGamma} \frac{\Gamma(n+2+2\RE w)}{(n+1)!} = N_{n+1}(-2\RE w -1) \Gamma(2 \RE w +1). \end{equation}

Since $n+2+ 2 \RE w >0$, we can use a table of integrals \cite[p. 30]{TableofIntegrals:1954} to show that \begin{align} \frac{1}{2\pi}\int_{-\infty}^{\infty}
 \frac{\left|\Gamma\left(\frac{n+2+2 \RE w}{2}+i\alpha\right)\right|^2}
 {\Gamma\left(n+2+2 \RE w\right)}e^{-i \ln\left(\frac{t}{s}\right)\alpha} 
 d \alpha  &= \left[\frac{1}{2}\sech \left(\frac{\ln\left(\frac{t}{s}\right)}{2}\right)\right]^{n+2+2 \RE w} \notag \\ &=  \left[ \frac{\sqrt{st}}{t+s}\right]^{n+2+2 \RE w}. \label{sech}
  \end{align}  
  By applying (\ref{SwitchingGamma}) and (\ref{sech}) to (\ref{ipsofar}) and recalling property (\ref{Newtonsum}) of the Newton polynomials, we obtain
 \begin{align*}  \left<T_{t^z}^*T_{s^z}k_w, k_w\right>_{H^2(\mu)} &= \sum_{n=-1}^{\infty} (st)^{\frac{n}{2}}\left(\frac{s}{t}\right)^{i \IM w}N_{n+1}(-2\RE w -1) \left[ \frac{\sqrt{st}}{t+s}\right]^{n+2+2 \RE w}  \displaybreak[1]\\
 &=  \frac{s^w t^{\overline{w}}}{(t+s)^{2 \RE w+1}} \sum_{n=-1}^{\infty} N_{n+1}(-2 \RE w -1)\left(\frac{st}{t+s}\right)^{n+1} \displaybreak[1] \\
 &  = \frac{s^w t^{\overline{w}}}{(t+s)^{2 \RE w +1}} \left(1-\frac{st}{t+s}\right)^{-2 \RE w -1} \\
 & = \frac{s^wt^{\overline{w}}}{(s+t-st)^{2 \RE w +1}}. \qedhere
 \end{align*}
\end{proof}

We can apply the previous lemma to determine a set of points that is contained in the joint approximate point spectrum of $\{T_{s_1^z}, \ldots, T_{s_n^z}\}.$

\begin{lemma}\label{whatwemusthave} Let $n \in \IN$ and  $0 < s_1, s_2, \ldots, s_n <1$ with $s_j \neq s_k$
if $j \neq k$.  Then $$\overline{\left\{\left(s_1^{-\frac{1}{2}+iy},
 \ldots, s_n^{-\frac{1}{2}+iy}\right) : y \in
\IR \right\}} \cup \left\{(0, \ldots, 0)\right\} \subseteq
\sigma_{ap}(T_{s_1^z}, \ldots, T_{s_n^z}).$$ \end{lemma}

\begin{proof}  To show that $(0, \ldots, 0) \in \sigma_{ap}(T_{s_1^z}, \ldots, T_{s_n^z})$, let $\omega_{\ell}=\frac{\ell}{2}$ for $\ell \in \IN$, and consider the sequence $\{k_{\omega_{\ell}}\}_{\ell \in \IN}$ of normalized reproducing kernel functions for $\mathcal{N}.$   Then by Lemma \ref{innerprodcalc}, \begin{equation*} \left|\left|T_{s_j^z}k_{\omega_{\ell}}\right| \right|^2_{H^2(\mu)} = \frac{s_j^{\ell}}{(2s_j-s_j^2)^{\ell + 1}}  =  \frac{1}{(2s_j-s_j^2)}\left(\frac{1}{2-s_j}\right)^{\ell} \end{equation*}  for $1 \leq j \leq n$. Thus, $||T_{s_j^z}k_{\omega_{\ell}}||_{H^2(\mu)} \rightarrow 0$ as $\ell \rightarrow \infty$ for all $j \in \{1, \ldots, n\}$, and $(0, \ldots, 0) \in \sigma_{ap}(T_{s_1^z}, \ldots, T_{s_n^z}).$ 

Now consider $\lambda=(\lambda_1, \ldots, \lambda_n) \in  \overline{\left\{\left(s_1^{-\frac{1}{2}+iy},
 \ldots, s_n^{-\frac{1}{2}+iy}\right) : y \in
\IR\right\}}$. Fix a sequence $\{y_\ell\}_{\ell=1}^{\infty}$ of real numbers such that $$\lim_{\ell \rightarrow \infty} 
\left(s_1^{-\frac{1}{2}+iy_{\ell}}, \ldots, s_n^{-\frac{1}{2}+iy_{\ell}}\right) = \lambda.$$
Notice that since  $\lim_{\ell \rightarrow \infty}s_j^{-\frac{1}{2} + iy_{\ell}} =\lambda_j$ for $1\leq j \leq n$, it is required that $|\lambda_j|=s_j^{-\frac{1}{2}}$ for all $j \in \{1, \ldots, n\}.$ 

Let $w_{\ell}=-\frac{1}{2}+\frac{1}{\ell}+iy_{\ell}$ for $\ell \in \IN,$ and consider the sequence $\{k_{w_\ell}\}_{\ell \in \IN}$ of normalized reproducing kernel functions for $\mathcal{N}$ . 
 For $1 \leq j \leq n$, \begin{align*} \left|\left| (T_{s_j^z}-\lambda_jI)k_{w_{\ell}}\right|\right|^2_{H^2(\mu)} 
 & = \frac{s_j^{w_{\ell}} s_{j}^{\overline{w_{\ell}}}}{(2s_j-s_j^2)^{2\RE w_{\ell}+1}} -2\RE \overline{\lambda_j}s^{w_{\ell}} + |\lambda_j|^2 \\ \displaybreak[0]
 & = \frac{1}{s_j}\left(\frac{1}{2-s_j}\right)^{\frac{2}{\ell}} -2 \text{Re}\overline{\lambda_j}s_j^{-\frac{1}{2}+iy_{\ell}}s_j^{\frac{1}{\ell}} +\frac{1}{s_j} \end{align*}
by Lemma \ref{innerprodcalc}.  
Since \begin{equation*} \lim_{\ell \rightarrow \infty}\left(\frac{1}{2-s_j}\right)^{\frac{2}{\ell}} =1 = \lim_{\ell \rightarrow \infty} s_j^{\frac{1}{\ell}}\end{equation*} and  \begin{equation*}\lim_{\ell \rightarrow \infty}\overline{\lambda_j}s_j^{-\frac{1}{2}+iy_{\ell}}=|\lambda_j|^2=s_j^{-1}, \end{equation*} we obtain that $ \left|\left| (T_{s_j^z}-\lambda_jI)k_{w_{\ell}}\right|\right|_{H^2(\mu)} \rightarrow  0$ as $\ell \rightarrow \infty$ for $1 \leq j \leq n$.  Hence, $\lambda \in \sigma_{ap}(T_{s_1^z}, \ldots, T_{s_n^z}).$
\end{proof}


We want to show the the points specified in Lemma \ref{whatwemusthave} are the only points in the joint approximate point spectrum.  We need an additional tool to help us exclude some of the other points.
 Given a measure space $(X, \nu)$, let
$\varphi_1, \varphi_2, \ldots, \varphi_n \in L^{\infty}(X, \nu)$.
The joint essential range of $\varphi_1,  \ldots,
\varphi_n$, denoted $\mathcal{E}_{\nu}(\varphi_1,  \ldots,
\varphi_n)$, is the set of all $(\lambda_1,  \ldots, \lambda_n) \in \IC^n$
such that for all $\varepsilon > 0$, \begin{equation*} \nu\left(\left\{z \in X:
\sum_{j=1}^{n} \left|\varphi_j(z)-\lambda_j\right| <
\varepsilon\right\}\right) >0.\end{equation*}  Note that, for $n=1$, $\mathcal{E}_{\nu}(\varphi_1)$ is the usual essential range of $\varphi_1.$  In our case, $X=\Omega$ and $\nu=\mu$, the measure defined by (\ref{measure}).

\begin{lemma}\label{DashApp1} Let $n \in \IN$ and  $0 < s_1, s_2, \ldots, s_n <1$ with $s_j \neq s_k$
if $j \neq k$.  Then $\sigma_{ap}\left(T_{s_1^z}, \ldots,
T_{s_n^z}\right) \subseteq \mathcal{E}_{\mu}\left(s_1^z,  \ldots, s_n^z\right).$
\end{lemma}
\begin{proof}  The arguments in this proof follow closely those of the proof of
Theorem 5.2 in \cite{Dash:1973}.  

Let $(\lambda_1,  \ldots, \lambda_n) \in \IC^n \setminus \mathcal{E}_{\mu}\left(s_1^z, \ldots, s_n^z\right)$.  We want to show that there exist $\psi_1, \psi_2, \ldots, \psi_n \in L^{\infty}(\mu)$ such that $\sum_{j=1}^n \psi_j(s_j^z-\lambda_j) = 1$ $\mu$-almost everywhere.  In that case, $\sum_{j=1}^{n} T_{\psi_j}(T_{s_j^z} -\lambda_jI) = T_{\sum_{j=1}^{n} \psi_j(s^z-\lambda_j)} = I$ and thus $(\lambda_1,  \ldots, \lambda_n) \notin \sigma_{ap}\left(T_{s_1^z}, \ldots, T_{s_n^z}\right).$

Suppose not, i.e.\ suppose that for all $\psi_1,  \ldots, \psi_n \in L^{\infty}(\mu)$, $\sum_{j=1}^{n} \psi_j(s_j^z-\lambda_j)$ is not invertible in $L^{\infty}(\mu)$.  This is equivalent to saying that, for all $\varepsilon >0$ and $\psi_1, \ldots, \psi_n \in L^{\infty}(\mu)$,  $$\mu\left(\left\{z \in \Omega :  \left|\sum_{j=1}^{n} \psi_j(z)\left(s_j^z-\lambda_j\right)\right| < \varepsilon\right\}\right) >0.$$ By setting $\psi_j=\overline{s_j^z-\lambda_j}$, we obtain that  \begin{equation*}\mu\left(E_{1}(\varepsilon):=\left\{z \in \Omega :   \sum_{j=1}^{n} |s_j^z-\lambda_j|^2 < \varepsilon\right\}\right)>0 \end{equation*}  for all $\varepsilon >0.$
But since $\lambda \notin \mathcal{E}_{\mu}(s_1^z, \ldots, s_n^z)$, there exists $\varepsilon^{\prime} >0$ such that $$\mu\left(E_2(\varepsilon^{\prime}):=\left\{z \in \Omega:  \sum_{j=1}^{n}\left|s_j^z-\lambda_j\right| < \varepsilon^{\prime}\right\}\right) =0.$$ 
 However, for all $z \in \IC$, \begin{equation*} \sum_{j=1}^n |s_j^z-\lambda_j| \leq \left(\sum_{j=1}^n |s_j^z-\lambda_j|^2\right)^{\frac{1}{2}}\left(\sum_{j=1}^{n} 1 \right)^{\frac{1}{2}}=\left(n\sum_{j=1}^n |s_j^z-\lambda_j|^2\right)^{\frac{1}{2}} \end{equation*} by the Schwarz inequality. 
  Thus, if $\varepsilon = \frac{(\varepsilon^{\prime})^2}{n}$, $E_1(\varepsilon) \subseteq E_2(\varepsilon^{\prime})$ and  $\mu(E_1(\varepsilon))=0$, which is a contradiction.  \end{proof}

  \begin{corollary}\label{essrangeincl1}  Let $0 < s < 1$. Then $$\sigma_{ap}(T_{s^z}) \subseteq \left\{s^w : \RE w= \frac{m}{2}, m = -1, 0, 1, \ldots\right\} \cup \{0\}.$$ \end{corollary}
  
  We now temporarily restrict to the case of one operator.  To be able to determine when an operator $T_{s^z}-\lambda I$ is bounded below, we need to obtain a certain upper bound on the norms of functions in $H^2(\mu)$.
  
  \begin{lemma}\label{boundingwoapiece} If $f \in H^2(\mu)$ and $m \in \IN \cup \{0\}$, then \begin{equation}\label{boundbelow} ||f||_{H^2(\mu)}^2 \leq (m+3)  \sum_{\substack{n=-1 \\ n\neq m}}^{\infty} \int_{-\infty}^{\infty} \left|f\left(\frac{n}{2} + iy\right)\right|^2 \frac{\left|\Gamma\left( \frac{n}{2}+1+iy\right)\right|^2 dy}{2\pi(n+1)!}. \end{equation} \end{lemma}

\begin{proof}

Let $f \in H^2(\mu)$ and $m \in \IN \cup \{0\}$.  By the relationship between $H^2(\mu)$ and $\mathcal{N}$, there exists $g\in \mathcal{N}$ with $f=g$ $\mu$-almost everywhere.  Then $F(z):=g(z)\Gamma(z+1) \in H^2\left(\frac{-1}{2}, \frac{m+1}{2}\right)$
 by Theorem \ref{MRRThm}, and hence
$F \in H^2\left(\frac{m-1}{2}, \frac{m+1}{2}\right)$. 
  By Lemma \ref{three lines},
   \begin{equation}\label{usingtl} \left|\left|F\left(\frac{m}{2}+  i y\right)\right|\right|_2 \leq  \left|\left|F\left(\frac{m-1}{2} +  i y\right)\right|\right|_2^{\frac{1}{2}}\cdot  \left|\left|F\left(\frac{m+1}{2} +  i y\right)\right|\right|_2^{\frac{1}{2}},\end{equation}  where $|| \cdot ||_2$ indicates the norm in $L^2(\IR)$ and $y$ is the variable.
By manipulating (\ref{usingtl}), we obtain
  \begin{align}  \frac{\left|\left|F\left(\frac{m}{2}+ i y\right)\right|\right|_2^2}{2 \pi (m+1)!} 
  &  \leq  \max \left\{ \frac{\left|\left|F\left(\frac{m-1}{2} + i y\right)\right|\right|_2^2}{(m+1)2\pi(m!)}, 
  \frac{(m+2) \left|\left|F\left(\frac{m+1}{2} +  i y\right)\right|\right|_2^2}{2\pi(m+2)!} \right\}  \notag \\
  &\leq  (m+2)\left( \frac{\left|\left|F\left(\frac{m-1}{2} + i y\right)\right|\right|_2^2}{2\pi(m!)}+
  \frac{\left|\left|F\left(\frac{m+1}{2} + i y\right)\right|\right|_2^2}{2\pi(m+2)!} \right). \label{goodbound}
 \end{align}
  By filling in the definition of $F$, (\ref{goodbound}) becomes 
   \begin{align}  &\int_{-\infty}^{\infty} \left|g\left(\frac{m}{2}+iy\right)\right|^2 \frac{\left|\Gamma\left(\frac{m}{2}+1 + iy\right)\right|^2 }{2\pi(m+1)!} dy \notag \\
  & \qquad \qquad \leq (m+2)\left(\int_{-\infty}^{\infty} \left|g\left(\frac{m-1}{2}+iy\right)\right|^2\frac{\left|\Gamma\left(\frac{m-1}{2}+1 +iy\right)\right|^2 }{2\pi(m!)} dy\right. \notag \\
  & \qquad \qquad  \qquad \qquad \quad +  \left. \int_{-\infty}^{\infty} \left|g\left(\frac{m+1}{2}+iy\right)\right|^2\frac{\left|\Gamma\left(\frac{m+1}{2}+1 +iy\right)\right|^2}{2\pi(m+2)!}dy\right). \label{realwork}   \end{align}  By replacing the expression in the parentheses on the right hand side of (\ref{realwork}) by \begin{equation}\label{RHS}\sum_{\substack{n=-1 \\ n \neq m}}^{\infty} \int_{-\infty}^{\infty} \left|g\left(\frac{n}{2} + iy\right)\right|^2 \left|\Gamma\left( \frac{n}{2}+1+iy\right)\right|^2 \frac{dy}{2\pi(n+1)!} \end{equation} and then adding (\ref{RHS}) to both sides, we obtain (\ref{boundbelow}) with $f$ replaced by $g$.  Since $f=g$ $\mu$-almost everywhere, we have (\ref{boundbelow}) as written in the statement of the theorem. \end{proof}

We can use the preceding results to identify the approximate point spectrum of a single operator of the form $T_{s^z}.$

  \begin{proposition}\label{oneoperatorspectra}  If $0 < s <1,$ then the approximate point spectrum of $T_{s^z}$ is $$\sigma_{ap}(T_{s^z})=\left\{s^{-\frac{1}{2}+iy} : y \in \IR\right\} \cup \{ 0\}= s^{-\frac{1}{2}}\IT \cup \{0\}.$$\end{proposition}
 
 \begin{proof} By Lemma \ref{whatwemusthave} and Corollary \ref{essrangeincl1}, \begin{align*}\left\{ s^{-\frac{1}{2}+iy} : y \in \IR \right\} \cup \{0\} & \subseteq \sigma_{ap}(T_{s^z}) \subseteq \left\{s^w : \RE w= \frac{m}{2}, m = -1, 0, 1, \ldots\right\} \cup \{0\}.\end{align*}
 
 Let $z_0=\frac{m}{2}+iy_0$ for $m \in \IN \cup \{0\}$ and $y_0 \in \IR$.  We want to show that $s^{z_0} \notin \sigma_{ap}(T_{s^z}).$
 Note that for $z=\frac{n}{2}+iy$, where $y \in \IR$ and $n \in \IN \cup \{0,-1\}$ with $n \neq m$, \begin{equation*} \left|s^z-s^{z_0}\right| \geq \left|s^{\frac{n}{2}}-s^{\frac{m}{2}}\right| 
 \geq s^{\frac{m}{2}} \cdot \min\left\{\left|s^{\frac{1}{2}} -1\right|, \left|s^{-\frac{1}{2}}-1\right|\right\} :=\lambda_{s,m}. \end{equation*}
  By Lemma \ref{boundingwoapiece} and the definition of $\mu$,  \begin{align*} \left|\left|(T_{s^z}-s^{z_0}I)f\right|\right|_{H^2(\mu)}^2 
 & = \int_{\Omega} \left|s^{z}-s^{z_0}\right|^2\left|f\left(z\right)\right|^2
d\mu \displaybreak[0] \\& \geq  \sum_{\substack{n=-1 \displaybreak[0]\\ n \neq m}}^{\infty}
   \int_{-\infty}^{\infty} \lambda_{s,m}^2\left|f\left(\frac{n}{2}+iy\right)\right|^2\frac{\left|\Gamma\left(\frac{n}{2}+1+iy\right)\right|^2dy}{2\pi(n+1)! } \\ \displaybreak[0]
    & \geq  \frac{\lambda_{s,m}^2}{m+3}||f||_{H^2(\mu)}^2 \end{align*} for all $f \in H^2(\mu).$  Thus, $T_{s^z}-s^{z_0}I$ is bounded below, so $s^{z_0} \notin \sigma_{ap}(T_{s^z})$. 
 \end{proof}
 
 We now return to the case of considering an arbitrary finite number of operators.  By Theorem \ref{oneoperatorspectra} and (\ref{inclusionindirectsum}), \begin{equation}\label{inclinsum2} \sigma_{ap}(T_{s_1^z}, \ldots, T_{s_n^z}) \subseteq \mathcal{E}_{\mu}(s_1^z,
\ldots, s_n^z) \cap \prod_{j=1}^n \left\{
 s_{j}^{-\frac{1}{2}} \IT\cup
\{0\}\right\}.\end{equation}  We study the space on the right hand side of (\ref{inclinsum2}) through the following lemmas:

\begin{lemma}\label{cannotmix}Let $n \in \IN$ and  $0 < s_1, s_2, \ldots, s_n <1$ with $s_j \neq s_k$
if $j \neq k$.  Then \begin{equation*}\mathcal{E}_{\mu}(s_1^z,
 \ldots, s_n^z) \cap \prod_{j=1}^n \left\{
s_j^{-\frac{1}{2}}\IT \cup \{0\}\right\} =\left\{\mathcal{E}_{\mu}(s_1^z,
 \ldots, s_n^z) \cap \prod_{j=1}^n 
s_j^{-\frac{1}{2}}\IT   \right\} \cup \{(0, \ldots,0)\}.\end{equation*}
 \end{lemma}

\begin{proof} The statement is trivially true if $n=1$, so we assume $n \geq 2.$   
 Suppose $(\lambda_1, \ldots, \lambda_n) \in  \mathcal{E}_{\mu}(s_1^z,
 \ldots, s_n^z) \cap \prod_{j=1}^n \left\{
s_j^{-\frac{1}{2}}\IT \cup \{0\}\right\}$ 
with $\lambda_{j_0}=0$ for some $j_0 \in \{1, \ldots, n\}$ and
 $\lambda_{j_1} \neq 0$ for some $j_1 \in \{1, \ldots, n\} \setminus \{j_0\}.$
  Since $\lambda_{j_1} \neq 0,$ there exists $y_{j_1} \in \IR$ such that $\lambda_{j_1}=s_{j_1}^{-\frac{1}{2}+iy_{j_1}}$.

Let $\varepsilon = \min\left\{1, s_{j_1}^{-\frac{1}{2}}-1\right\}.$  Notice that if $z \in \IC$ with $\text{Re} \, z = \frac{m}{2}$ for some $m \in \{ -1, 0, 1, \ldots\}$, then, for any $j$,\begin{equation}\label{distance1} \left|s_j^z-\lambda_{j_1}\right| = \left|s_j^z-s_j^{-\frac{1}{2}+iy_{j_1}}\right| \geq s_j^{-\frac{1}{2}}-|s_j^z| \geq s_j^{-\frac{1}{2}}-1 \geq \varepsilon  \end{equation} if $  m \geq 0$, while
\begin{equation}\label{distance2}|s_{j_0}^z-0| \geq 1 \geq \varepsilon  \end{equation} if $m = -1.$
   Therefore, \begin{equation*} \mu\left(\left\{ z\in \Omega: \sum_{j=1}^n |s_j^z - \lambda_j| < \varepsilon\right\}\right) =0,\end{equation*} 
which contradicts the fact that $(\lambda_1, \ldots, \lambda_n) \in \mathcal{E}_{\mu}(s_1^z,  \ldots, s_n^z)$ and thus proves the lemma. \end{proof}
    
    \begin{lemma}\label{lastjointpiece}Let $n \in \IN$ and  $0 < s_1, s_2, \ldots, s_n <1$ with $s_j \neq s_k$
if $j \neq k$.  Then \begin{equation*}\mathcal{E}_{\mu}\left(s_1^z, \ldots, s_n^z\right) \cap \prod_{j=1}^n s_j^{-\frac{1}{2}} \IT \subseteq \overline{\left\{\left(s_1^{-\frac{1}{2}+iy}, \ldots, s_n^{-\frac{1}{2}+iy} \right) : y \in \IR \right\}}.\end{equation*} \end{lemma}  

\begin{proof}
For $j=1, 2, \ldots, n$, consider the functions $\psi_{s_j}: \IR \rightarrow \IC$ defined by $\psi_{s_j}(y)=s_j^{-\frac{1}{2}+iy}$.  
    We can view $\psi_{s_j}$ as $s^z |_{\text{Re} \, z = -\frac{1}{2}}$.  
    By restricting the measure $\mu$ to the line $\RE z= -\frac{1}{2}$ in a similar way, we obtain a measure $\hat{\mu}$ on $\IR$ given by $$d\hat{\mu} = \frac{\left|\Gamma\left(\frac{1}{2}+iy\right)\right|^2}{2 \pi} dy.$$
     Notice that $\psi_{s_j} \in L^{\infty}(\IR, \hat{\mu})$ for all $1 \leq j \leq n$.
     Using the idea behind equation (\ref{distance1}) in this setting, it is easy to show that $$\mathcal{E}_{\mu}\left(s_1^z, \ldots, s_n^z\right) \cap \prod_{j=1}^n s_j^{-\frac{1}{2}}\IT \subseteq \mathcal{E}_{\hat{\mu}}\left(\psi_{s_1}, \ldots \psi_{s_n}\right).$$  Since $\hat{\mu}$ is mutually absolutely continuous with respect to Lebesgue measure and $\psi_{s_j}$ is continuous for $1 \leq j  \leq n$, $$\mathcal{E}_{\hat{\mu}}\left(\psi_{s_1}, \ldots \psi_{s_n}\right) = \overline{\left\{\left(s_1^{-\frac{1}{2}+iy}, \ldots, s_n^{-\frac{1}{2}+iy} \right) : y \in \IR \right\}},$$ which proves the lemma.\end{proof}

    Combining the results of Lemmas \ref{whatwemusthave}, \ref{DashApp1}, \ref{cannotmix}, and \ref{lastjointpiece}, we obtain the following theorem that includes the results of Proposition \ref{oneoperatorspectra} as a special case.
 \begin{theorem}\label{JAP1} Let $n \in \IN$ and  $0 < s_1, s_2, \ldots, s_n <1$ with $s_j \neq s_k$
if $j \neq k$.  Then \begin{equation*}\sigma_{ap}(T_{s_1^z},  \ldots, T_{s_n^z}) = \{(0,\ldots, 0)\} \cup \overline{\left\{\left(s_1^{-\frac{1}{2}+iy}, \ldots, s_n^{-\frac{1}{2}+iy}\right) : y \in \IR \right\}}.\end{equation*} \end{theorem}

The following corollary of Theorem \ref{JAP1} is immediate from Theorem \ref{Buncefinite}.

 \begin{corollary} Let $n \in \IN$ and  $0 < s_1, s_2, \ldots, s_n <1$ with $s_j \neq s_k$
if $j \neq k$.  Let $\mathcal{C}_{\{s_1, \ldots, s_n\}}$ be the commutator ideal of $C^*\left(T_{s_1^z},  \ldots, T_{s_n^z}\right).$  Then \begin{equation*}   
C^*\left(T_{s_1^z}, \ldots, T_{s_n^z}\right)  / \mathcal{C}_{\{s_1, \ldots, s_n\}}  \cong  C\left(\overline{\left\{\left(s_1^{-\frac{1}{2}+iy},  \ldots, s_n^{-\frac{1}{2}+iy}\right) : y \in \IR \right\}}\right) \oplus \IC. \end{equation*} \end{corollary}
 
 \subsection{Investigating the Shapes of the Sets $\sigma_{ap}(T_{s_1^z}, \ldots, T_{s_n^z})$}
 
 We would like to use Theorem \ref{Bunceinfinite} to determine the structure of $C^*(\{T_{s^z} : 0 < s < 1 \})$ modulo its commutator ideal.  Thus, we need to understand the shapes of the joint approximate point spectra of all finite subsets of $\{T_{s^z} : 0 < s < 1\}$ so that we can compute the needed inverse limit.  The structures of these spectra depend on the relations between the numbers $\ln(s_1), \ldots, \ln (s_n)$ via Kronecker's Theorem. The following version of Kronecker's Theorem is included in \cite{LevitanZhikov:1982}.

\begin{theorem}[Kronecker's Theorem] Let $\alpha_1, \alpha_2, \ldots, \alpha_n, \theta_1, \theta_2, \ldots, \theta_n$ be arbitrary real numbers.  For the system of inequalities \begin{equation*} |\alpha_k t-\theta_k| < \delta \mod 2\pi \, \, \, (k=1, 2, \ldots, n) \end{equation*} to have consistent real solutions for any arbitrarily small positive number $\delta$, it is necessary and sufficient that every time the relation $k_1\alpha_1+ k_2\alpha_2+ \ldots +k_n\alpha_n =0$ holds, where $k_1, k_2, \ldots, k_n$ are integers, we have the congruence \begin{equation*} k_1\theta_1 + k_2\theta_2 + \ldots k_n\theta_n \equiv 0 \mod 2\pi. \end{equation*} \end{theorem}

The simplest case of Kronecker's Theorem is when $\alpha_1, \alpha_2, \ldots, \alpha_n$ are linearly independent over $\IZ.$  
A finite collection $\{\alpha_1, \alpha_2, \ldots, \alpha_n\}$ of real numbers is linearly independent over $\IZ$ if $k_1 \alpha_1 + k_2\alpha_2 + \ldots +k_n\alpha_n=0$ with $k_1, k_2, \ldots, k_n \in \IZ$ if and only if $k_1=k_2=\ldots =k_n=0$.     In this case, we can combine Theorem \ref{JAP1} with a method used by 
B\"ottcher, Karlovich, and Spitkovsky in \cite[Corollary 1.13]{BottcherKarlovichSpitkovsky:2002} to straight-forwardly show the following result:

\begin{corollary}\label{japli}  Let $n \in \IN$ and $0 < s_1, s_2, \ldots, s_n <1$ with $s_j \neq s_k$ if $j \neq k$.  If the numbers $\ln(s_1), \ln(s_2), \ldots, \ln(s_n)$ are linearly independent over $\IZ$, then  \begin{equation*}\sigma_{ap}(T_{s_1^z}, \ldots, T_{s_n^z})= \{(0, \ldots, 0)\} \cup \left(s_1^{-\frac{1}{2}}\IT \times s_2^{-\frac{1}{2}} \IT \times \ldots \times s_n^{-\frac{1}{2}}\IT\right). \end{equation*} \end{corollary}

The other case in which the shape of $\sigma_{ap}\left(T_{s_1^z},  \ldots, T_{s_n^z}\right)$ is simple to determine is the case where all of the $\ln(s_j)$ are rational multiples of each other.

\begin{lemma}\label{japrm} Let $n \in \IN$ and $0 < s_1, s_2, \ldots, s_n <1$ with $s_j \neq s_k$ if $j \neq k$.  Suppose there exist  $a_2, \ldots, a_n, b_2, \ldots, b_n$ such that $\ln(s_j)=\frac{a_j}{b_j}\ln(s_1)$ and $\gcd(a_j, b_j)=1$  for $j=2, \ldots, n.$  Then the range of $\left(s_1^{-\frac{1}{2} +iy},  \ldots, s_n^{-\frac{1}{2}+ iy}\right)$ is closed, so \begin{equation*}\sigma_{ap}(T_{s_1^z},  \ldots, T_{s_n^z}) = \{(0,\ldots, 0)\} \cup \left\{\left(s_1^{-\frac{1}{2} +iy}, \ldots, s_n^{-\frac{1}{2}+ iy}\right) : y \in \IR \right\}.\end{equation*}
\end{lemma}

\begin{proof} Let $M=\text{lcm}[b_2, \ldots, b_n]$.  Let $\psi: \IR \rightarrow \IC$ be defined by $$\psi(y)=\left(s_1^{-\frac{1}{2} +iy}, \ldots, s_n^{-\frac{1}{2}+ iy}\right).$$
Clearly, $\psi$ is a continuous function. 
 We want to show that $\psi$ is periodic with period $\frac{2M\pi}{-\ln(s_1)}. $
  Let $y \in \IR$ and $k \in \IZ$.  
  Then it is clear that \begin{equation*} s_1^{-\frac{1}{2} +i\left(y+k\left(\frac{2M\pi}{-\ln(s_1)}\right)\right)} = s_1^{-\frac{1}{2}+iy}e^{\ln(s_1)ik\left(\frac{2M\pi}{-\ln(s_1)}\right)} = s_1^{-\frac{1}{2}+iy}. \end{equation*} 
  If $2 \leq j \leq n$, then \begin{equation*} s_j^{-\frac{1}{2}+i\left(y+k\left(\frac{2M\pi}{-\ln(s_1)}\right)\right)}= s_j^{-\frac{1}{2}+iy}e^{\frac{a_j}{b_j}\ln(s_1) i k\left(\frac{2M\pi}{-\ln(s_1)}\right)}  = s_j^{-\frac{1}{2}+iy} \end{equation*} since $b_j$ divides $M$.  Thus, $\psi$ is a periodic function with period $\frac{2M\pi}{-\ln(s_1)}$ and $\psi(\IR)$ $ =\psi\left(\left[0,  \frac{2M\pi}{-\ln(s_1)}\right]\right)$, which is a closed set.
  \end{proof}
 
  Notice that in the case of Lemma \ref{japrm}, $\sigma_{ap}\left(T_{s_1^z},  \ldots, T_{s_n^z}\right)$ consists of the point $(0, \ldots, 0)$ 
  and a closed curve that is homeomorphic to $\IT.$  Thus, the two distinct cases described in Corollary \ref{japli} and Lemma \ref{japrm} lead to joint approximate point spectra that are not homeomorphic to each other. 
  
 The preceding investigations fully determine the  possible shapes of the joint approximate point spectrum of two operators from $\{T_{s^z} : 0 < s < 1\}$.  To consider three or more operators, one must investigate a larger number of possible relationships between $\ln(s_1), \ln(s_2), \ldots, \ln(s_n),$ and the determination of the joint approximate spectrum requires the full version of Kronecker's Theorem. Even for three operators, the calculations quickly become quite complicated.

\subsection{Calculating $\sigma_{ap}(W_{\chi_{\alpha_1}}
, \ldots, W_{\chi_{\alpha_n}})$}

  Due to the complexity of the joint approximate point spectra, it appears impractical to determine directly the structure of the inverse limit of all joint approximate point spectra of finite subsets of $\{T_{s^z}: 0 < s < 1\}.$  Instead, we will build a similar framework for Toeplitz operators on $H^2(\IR)$ and determine the desired inverse limit by connecting the spectral results for the two collections of operators and recalling that the structure of $C^*(W_{\chi_{\alpha}})$, modulo its commutator ideal, is described by Theorem \ref{CDiso}.
  
 For these purposes, we want to understand the joint approximate spectra of finite subsets of $\{W_{\chi_{\alpha}} : \alpha \in \IR^+\}.$  We begin with a result of Dash that identifies the joint approximate point spectrum of any finite collection of multiplication operators on $L^2(\IR).$

\begin{theorem}\label{Dash2}{\rm \cite[Theorems 5.2, 5.3]{Dash:1973}} If $n \in \IN$ and $\varphi_1, \varphi_2, \ldots, \varphi_n \in L^{\infty}(\IR)$, then \begin{equation*}\sigma_{ap}\left(M_{\varphi_1},  \ldots, M_{\varphi_n}\right) = \mathcal{E}_{m}(\varphi_1, \ldots,  \varphi_n). \end{equation*} \end{theorem}

To use Theorem \ref{Dash2} in our setting, we apply the methods used in the proof of Theorem 1 in \cite{CoburnDouglas:1971} to multiple operators
 to prove the following lemma.  In the statement of the lemma, we restrict to maps $\varphi_1, \ldots, \varphi_n \in L^{\infty}(\IR) \cap H^2(\IR)$ since the joint approximate point spectrum is only defined for collections of commuting operators.  We could have alternately considered $\varphi_1, \ldots, \varphi_n \in L^{\infty} (\IR)$ with $\overline{\varphi_1}, \ldots, \overline{\varphi_n} \in H^2(\IR).$

\begin{lemma}\label{multinToep} Let $n \in \IN$ and  $\varphi_1, \varphi_2, \ldots, \varphi_n \in H^2(\IR) \cap L^{\infty}(\IR).$  Then $$\sigma_{ap}\left(M_{\varphi_1},  \ldots, M_{\varphi_n}\right) \subseteq \sigma_{ap}\left(W_{\varphi_1}, \ldots, W_{\varphi_n}\right).$$  \end{lemma}
\begin{proof}  
Note that it suffices to show that if $(0, \ldots, 0) \in \sigma_{ap}\left(M_{\varphi_1},  \ldots, M_{\varphi_n}\right)$, then $(0, \ldots, 0) \in \sigma_{ap}\left(W_{\varphi_1},  \ldots, W_{\varphi_n}\right).$
 For $1 \leq j \leq n$, consider the net of operators $\left\{B_{j,\alpha}\right\}_{\alpha \in \IR^+}$ in $\mathcal{B}(L^2(\IR))$ defined by
 \begin{eqnarray*} B_{j, \alpha}& = &M_{\chi_{\alpha}}^{*}W_{\varphi_j}P_mM_{\chi_{\alpha}}\\
 & = & M_{\chi_{\alpha}}^{*}P_mM_{\chi_{\alpha}}M_{\chi_{\alpha}}^{*}M_{\varphi_j}P_mM_{\chi_{\alpha}}\\  
 &=& (M_{\chi_{\alpha}}^{*}P_mM_{\chi_{\alpha}})M_{\varphi_j}(M_{\chi_{\alpha}}^{*}P_mM_{\chi_{\alpha}}). \end{eqnarray*}  
 In \cite{CoburnDouglas:1971}, Coburn and Douglas showed that, for all $j \in \{1, \ldots, n\}$, the net
 $\left\{B_{j,\alpha}\right\}_{\alpha \in \IR^+}$ converges to $M_{\varphi_j}$ in the strong operator topology on $\mathcal{B}(L^2(\IR)).$  They also showed that $\left\{M_{\chi_{\alpha}}^*P_m M_{\chi_{\alpha}}\right\}_{\alpha \in \IR^+}$ converges strongly to $I$, the identity operator on $L^{2}(\IR).$ 

Suppose $(0, \ldots, 0) \in \sigma_{ap}(M_{\varphi_1}, \ldots, M_{\varphi_n})$.  Let $\varepsilon >0$ be given.  Then there exists a unit vector $f \in L^2(\IR)$ such that \begin{equation*} ||M_{\varphi_j}f||_{2} < \frac{\varepsilon}{4}\end{equation*} for all $ 1 \leq j  \leq n$.  
Since $\left\{B_{j,\alpha}\right\}_{\alpha \in \IR^+}$ converges strongly to $M_{\varphi_j}$, there exists, for each $1 \leq j \leq n$, a number $\alpha_{j,0} \in \IR^+$ such that if $\alpha \geq \alpha_{j,0}$,
then \begin{equation}\label{newsmallvector} \left|\left|W_{\varphi_j}P_mM_{\chi_{\alpha}}f \right| \right|_2= \left|\left|M_{\chi_{\alpha}}^{*}W_{\varphi_j}P_mM_{\chi_{\alpha}} f \right|\right|_2 = ||B_{j, \alpha}f||_2 < \frac{\varepsilon}{2}  . \end{equation}  
We set $ \alpha_0 := \max_{1 \leq j \leq n} \left\{\alpha_{j,0}\right\}.$  Similarly, since the net $\left\{M_{\chi_{\alpha}}^*P_m M_{\chi_{\alpha}}\right\}_{\alpha \in \IR^+}$ converges strongly to $I$ and $f$ is a unit vector, there exists $\beta_0 \in \IR^+$, such that if $\alpha \geq \beta_0$, then \begin{equation}\label{othersmallvector} \left|\left|P_mM_{\chi_{\alpha}} f \right| \right|_2 = \left|\left|M_{\chi_{\alpha}}^*P_mM_{\chi_{\alpha}} f \right| \right|_2 > \frac{1}{2}.\end{equation}

Combining (\ref{newsmallvector}) and (\ref{othersmallvector}), we find that  if $\gamma \geq \alpha_0 + \beta_0$, then\begin{equation*} \left| \left|W_{\varphi_j} P_m M_{\chi_{\alpha}}f \right| \right|_2 < \frac{\varepsilon}{2} < \varepsilon \left|\left| P_m M_{\chi_{\alpha}} f\right| \right|_2. \end{equation*}  Thus, we can construct a sequence $\{f_{\ell}\}_{\ell=1}^{\infty}$ of unit vectors in $H^{2}(\IR)$ such that, for $1 \leq j \leq n,$  $||W_{\varphi_j}f_{\ell}||_{H^2(\mu)} \rightarrow 0$ as $\ell \rightarrow \infty$, so $(0, \ldots, 0) \in \sigma_{ap}(W_{\varphi_1}, \ldots, W_{\varphi_n}).$
\end{proof}

We now restrict our attention to the Toeplitz operators on $H^2(\IR)$ induced by the functions in $\{\chi_{\alpha} : \alpha \in \IR^+\}$.

\begin{theorem}\label{JAP2} Let $n \in \IN$ and $\alpha_1, \alpha_2, \ldots, \alpha_n \in \IR^+$ with $\alpha_j \neq \alpha_k$ if $j \neq k$.  Then \begin{equation*}
\sigma_{ap}\left(W_{\chi_{\alpha_1}}, \ldots, W_{\chi_{\alpha_n}}\right) = \overline{\left\{\left(\chi_{\alpha_1}(y), , \ldots, \chi_{\alpha_n}(y)\right): y \in \IR\right\}}. \end{equation*} \end{theorem}

\begin{proof}
By Theorem \ref{Dash2} and Lemma \ref{multinToep}, $$\mathcal{E}_{m}(\chi_{\alpha_1}, \ldots, \chi_{\alpha_n}) \subseteq \sigma_{ap}(W_{\chi_{\alpha_1}}, \ldots, W_{\chi_{\alpha_n}}).$$  The reverse inclusion can be proved by repeating the arguments of the proof of 
Lemma \ref{DashApp1}.  Since the functions $\chi_{\alpha_1}, \ldots, \chi_{\alpha_n}$ are continuous on $\IR$, \begin{equation*}\mathcal{E}_{m}(\chi_{\alpha_1}, \ldots, \chi_{\alpha_n}) =
 \overline{\left\{\left(\chi_{\alpha_1}(y),\ldots, \chi_{\alpha_n}(y)\right): y \in \IR\right\}}. \qedhere\end{equation*} 
 \end{proof}
      
      \section{\textsc{The Structure of} $C^*(\{C_{\varphi_s} : 0 < s <1 \})$ \textsc{Modulo the Commutator Ideal}}
      
     We will use the results from the previous section to show a connection between the set $\left\{T_{s_1^z}, \ldots, T_{s_n^z}\right\}$ in $\mathcal{B}(H^2(\mu))$ and a collection $\left\{W_{\chi_{\alpha_1}}, \ldots, W_{\chi_{\alpha_n}}\right\}$ in $\mathcal{B}(H^2(\IR)).$  This relationship will be a key ingredient in identifying
  the structure of $C^*(\left\{C_{\varphi_s}: 0 < s< 1\right\})$ modulo the commutator ideal.

 \begin{theorem} Let $n \in \IN$ and $\alpha_1, \alpha_2, \ldots, \alpha_n \in \IR^+$ with $\alpha_j \neq \alpha_k$ if $j \neq k$. Then \begin{equation*} \sigma_{ap}\left(e^{\frac{-\alpha_1}{2}}T_{e^{-\alpha_1 z}},  \ldots, e^{\frac{-\alpha_n}{2}}T_{e^{-\alpha_nz}}\right)=\{(0,\ldots, 0)\} \cup  \sigma_{ap}\left(W_{\chi_{\alpha_1}},  \ldots, W_{\chi_{\alpha_n}}\right). \end{equation*} \end{theorem}
 
\begin{proof}
By Theorem \ref{JAP1}, $\sigma_{ap}\left(T_{e^{-\alpha_1 z}}, \ldots, T_{e^{-\alpha_nz}}\right)$ equals 
 \begin{equation*}\left\{(0, \ldots, 0)\right\} \cup \overline{\left\{\left(e^{-\alpha_1(-\frac{1}{2}+iy)}, \ldots, e^{-\alpha_n(-\frac{1}{2}+iy)}\right) : y \in \IR \right\}}.\end{equation*}  
Then by simple arguments,
\begin{align*}  \sigma_{ap}\left(e^{\frac{-\alpha_1}{2}}T_{e^{-\alpha_1 z}},  \ldots, e^{\frac{-\alpha_n}{2}}T_{e^{-\alpha_nz}}\right) & = \{(0, \ldots, 0)\} \cup \overline{\left\{\left(e^{-\alpha_1 i y}, \ldots, e^{-\alpha_n iy}\right) : y \in \IR\right\}}\\ & = \{(0, \ldots, 0)\} \cup \overline{\left\{\left(\chi_{\alpha_1}(y), \ldots, \chi_{\alpha_n}(y)\right) : y \in \IR \right\}}. \end{align*} 
The result follows from Theorem \ref{JAP2}. \end{proof}
 
 We will  construct two inverse limits systems.  Consider the set $\mathcal{P}_{fin}^+$ of all non-empty finite subsets of $\IR^+$, which is partially ordered by inclusion.  For clarity of notation, if $F \in \mathcal{P}_{fin}^+$ has $n$ elements, then we write $F=\{\alpha_1, \alpha_2, \ldots, \alpha_n\}$, where $\alpha_1 < \alpha_2 < \ldots < \alpha_n.$ 
 
  For $F=\left\{\alpha_1, \alpha_2, \ldots, \alpha_n\right\} \in \mathcal{P}_{fin}^+$,we define two sets: \begin{equation*} \sigma_{F} = \sigma_{ap}\left(W_{\chi_{\alpha_1}},  \ldots, W_{\chi_{\alpha_n}}\right) \end{equation*} and \begin{equation*}\tilde{\sigma}_F= \sigma_{ap}\left(e^{-\frac{\alpha_1}{2}}T_{e^{-\alpha_1z}}, \ldots, e^{-\frac{\alpha_n}{2}}T_{e^{-\alpha_n z}}\right). \end{equation*}  Notice that, for all $F \in \mathcal{P}_{fin}^{+}$, $\tilde{\sigma}_F=\left\{(0,\ldots,0)\right\} \cup \sigma_F.$   We consider these sets as topological spaces in the relative topology from $\IC^{|F|}$.  We first investigate the connections between the open sets in the spaces. 
  
  \begin{lemma}\label{opensets}  Let $F \in \mathcal{P}_{fin}^+$.  Then the collections of open subsets of $\sigma_F$ and $\tilde{\sigma}_F$ are related in the following way: \begin{align*} & \left\{W : W \, \text{is open in} \,   \tilde{\sigma}_F\right\}\\ & \qquad \quad  =\left\{E : E \, \text{is open in } \, \sigma_{F}\right\} \cup \left\{ E \cup \{(0, \ldots, 0)\} : E \, \text{is open in } \, \sigma_{F}\right\}. \end{align*} \end{lemma}
  
  \begin{proof} Let $W$ be an open set in $\tilde{\sigma}_F$.  Then $W=V \cap \tilde{\sigma}_F$, where $V$ is an open subset of $\IC^{|F|}$, so $V \cap \sigma_{F}$ is open in $\sigma_{F}$.  But $V \cap \tilde{\sigma}_F = (V \cap \sigma_F) \cup (V \cap \{(0, \ldots, 0)\}.$  Thus, $W = V \cap \sigma_{F}$ if $(0, \ldots, 0) \notin W$, and $W= (V \cap \sigma_F) \cup \{(0, \ldots, 0)\}$ if $(0, \ldots, 0) \in W.$
  
 Let $E$ be an open set in $\sigma_F$.  Then $E=U \cap \sigma_F$ for some open set $U$ in $\IC^{|F|}.$  Let $\Delta_{\frac{1}{4}}(0)$ be the open polydisk in $\IC^{|F|}$ centered at $0$ with radius $\frac{1}{4}$.
Then $E =U \cap \left(\IC \setminus \overline{\Delta_{\frac{1}{4}} (0)}\right) \cap \sigma_F = U \cap \left(\IC \setminus\overline{\Delta_{\frac{1}{4}} (0)}\right) \cap \tilde{\sigma}_F,$ since all components of elements in $\sigma_F$ have modulus 1, so $E$ is open in $\tilde{\sigma}_F$. 
  Also, $E \cup \{(0, \ldots, 0)\} = (U \cup \Delta_{\frac{1}{4}}(0)) \cap \tilde{\sigma}_F$, so $E  \cup \{(0,\ldots, 0)\}$ is open in $\tilde{\sigma}_F$ as well.
 \end{proof}
  For $F, G \in \mathcal{P}_{fin}^+$ with  $F \subseteq G$, we define the map $\pi_{FG} : \sigma_G \rightarrow \sigma_F$ as the projection onto the coordinates coming from the elements of $F$.  This map is well-defined and surjective due to the properties of the joint approximate point spectrum.  The map $\tilde{\pi}_{FG}: \tilde{\sigma}_{G} \rightarrow \tilde{\sigma}_F$ is defined equivalently.  
 Note that if $F \subseteq G$ and $z \in \tilde{\sigma}_{G}$, then $\tilde{\pi}_{FG}(z)=\pi_{FG}(z)$ for $z \in \sigma_G$ and $\tilde{\pi}_{FG}(z)= (0, \ldots, 0) \in \IC^{|F|}$ for $z = (0, \ldots, 0) \in \IC^{|G|}.$
 
The maps $\pi_{FG}$ and $\tilde{\pi}_{FG}$ are clearly continuous. If $F \in \mathcal{P}_{fin}^+,$ then $\pi_{FF}$ and $\tilde{\pi}_{FF}$ are the identity transformations.  Also if $F, G$, and $H$ are in $\mathcal{P}_{fin}^+$ with $F \subseteq G \subseteq H$, then $\pi_{FG}\pi_{GH}=\pi_{FH}$ and $\tilde{\pi}_{FG}\tilde{\pi}_{GH}=\tilde{\pi}_{FH}.$  Thus, $\left\{\left\{\sigma_F\right\}, \left\{\pi_{FG}\right\}\right\}$ and $\left\{\left\{\tilde{\sigma}_F\right\}, \left\{\tilde{\pi}_{FG}\right\}\right\}$ are both inverse limit systems of topological spaces over $\mathcal{P}_{fin}^+.$
 
 We can construct the inverse limits of the these systems as subspaces of appropriate product spaces.    Recall that if we have a collection of topological spaces $\{X_{\alpha}\}_{\alpha \in \Lambda}$, then the product space $\prod_{\Lambda}X_{\alpha}$
   is the space of all functions $f: \Lambda \rightarrow \bigcup_{a \in \Lambda} X_{\alpha}$ that satisfy $f(\alpha) \in X_{\alpha}$ for all $\alpha \in \Lambda.$  For ease of notation, we will write elements of $\prod_{\Lambda}X_{\alpha}$ as nets $\left\{x_\alpha\right\},$ where $x_{\alpha}=f(\alpha)$ for all $\alpha \in \Lambda,$ instead of functions. A basis for the product topology on $\prod_{\Lambda} X_{\alpha}$ is \begin{equation*} \left\{\prod_{\Lambda} U_{\alpha} : \begin{array}{c}U_{\alpha} \, \, \text{is open in} \, \, X_{\alpha} \, \, \text{for all} \, \, \alpha \in \Lambda, \\ \text{all but finitely many of the}\, \,  U_{\alpha} \, \, \text{are equal to} \, \, X_{\alpha}\end{array}\right\}.  \end{equation*}  
   
   We construct the spaces \begin{equation*} a(\mathcal{W}) := 
\left\{ \left\{x_F\right\} \in \prod_{P_{fin}^+} \sigma_F : \pi_{FG}(x_G)=x_F \, \, \forall F \subseteq G\right\}   
\end{equation*} and
  \begin{equation*} a(\mathcal{T})
  := \left\{ \left\{\tilde{x}_F\right\} \in \prod_{\mathcal{P}_{fin}^+} \tilde{\sigma}_F : \tilde{\pi}_{FG}(x_G)=x_F \, \, \forall F \subseteq G\right\}. 
 \end{equation*}
Then, by standard facts about inverse limits systems of topological spaces, we have that $\displaystyle{a(\mathcal{W}) = \lim_{\leftarrow} \sigma_F}$ and $\displaystyle{a(\mathcal{T})= \lim_{\leftarrow} \tilde{\sigma}_F,}$ where $\displaystyle{\lim_{\leftarrow}}$ indicates the inverse limit taken over $\mathcal{P}_{fin}^+.$  The spaces $a(\mathcal{T})$ and $a(\mathcal{W})$ are non-empty, compact Hausdorff spaces in the relative topologies from  the product topologies on $\prod_{\mathcal{P}_{fin}^+} \sigma_F$ and $\prod_{\mathcal{P}_{fin}^+ } \tilde{\sigma}_F$, respectively.  For these and other facts about inverse limit systems of topological spaces, see \cite{HockingYoung:1961}.
 
 Notice that the net $\{\tilde{x}_F\}$, where $\tilde{x}_F=(0,\ldots, 0) \in \IC^{|F|}$ for all $F \in \mathcal{P}_{fin}^+$, is an element of $a(\mathcal{T}).$  For clarity, we denote this element by $\left\{(0,\ldots, 0)_F\right\}$ to distinguish it from the singleton set $\{(0,\ldots, 0)\}$ contained in an individual  $\tilde{\sigma}_F.$

 \begin{theorem} Let $a(\mathcal{W})$ and $a(\mathcal{T})$ be defined as above. Then \begin{equation}\label{equalplus0} a(\mathcal{T})=\left\{(0, \ldots, 0)_F\right\} \cup a(\mathcal{W}).\end{equation}   The open sets in $a(\mathcal{T})$ are the sets of the form $V$ and $\left\{(0,\ldots, 0)_F\right\} \cup V,$ where $V$ is an open set in $a(\mathcal{W}).$\end{theorem}
 
 \begin{proof} Let $\left\{x_F\right\} \in a(\mathcal{W})$.
  Then $x_F \in \sigma_F \subset \tilde{\sigma}_F$ for all  $F \in \mathcal{P}_{fin}^+$.
    If $F \subseteq G$, then $\tilde{\pi}_{FG}(x_G) = \pi_{FG}(x_G) = x_F$.  
    So $\{x_F\} \in a(\mathcal{T}).$  Also, as noted above, $\left\{(0, \ldots, 0)_F\right\} \in a(\mathcal{T})$.
 
 Let $\{\tilde{x}_F\} \in a(\mathcal{T}).$  Suppose $\tilde{x}_F \neq (0, \ldots, 0)$ for all $F \in \mathcal{P}_{fin}^+.$ 
 Then $\tilde{x}_F  
 \in \sigma_F$ for all $F \in \mathcal{P}_{fin}^+$.  Also if $F \subseteq G$, then $\pi_{FG}(\tilde{x}_G)=\tilde{\pi}_{FG}(\tilde{x}_G)=\tilde{x}_F$.  Thus $\{\tilde{x}_F\} \in a(\mathcal{W}).$  
 Now suppose there exists a set $F \in \mathcal{P}_{fin}^+$ such that $\tilde{x}_F=(0, \ldots, 0) $.  Since $\{\tilde{x}_F\} \in a(\mathcal{T})$ and, for all $G \in \mathcal{P}_{fin}^+$, all elements of $\tilde{\sigma}_G$ either have all components being zero or all non-zero components, it follows that $\tilde{x}_F=(0,\ldots, 0)$ for all $F \in \mathcal{P}_{fin}^+$ by the definitions of the maps $\tilde{\pi}_{FG}$.  Hence $\{\tilde{x}_F\} = \{(0, \ldots, 0)_F\}.$  Thus, we have shown (\ref{equalplus0}).  The relationship between the open sets of the two spaces is easy to show by using the bases for the topologies on the spaces and Lemma \ref{opensets}.   \end{proof}

  \begin{corollary}\label{contiso} Let $a(\mathcal{T})$ and $a(\mathcal{W})$ be defined as above. Then the map $\psi: C(a(\mathcal{T})) \rightarrow C(a(\mathcal{W})) \oplus \IC$, defined by $\psi(f)=\left(f|_{a(\mathcal{W})}, f(\{(0, \ldots 0)_F\})\right)$, is an isometric $*$-isomorphism. \end{corollary}
  
  \begin{proof} Since $a(\mathcal{W})$ can be viewed as a subspace of $a(\mathcal{T})$ with the relative topology, it is clear that $\psi$ is well-defined.  
  It is simple to show that $\psi$ is linear, multiplicative, $*$-preserving, and isometric. 
   We just need to check that $\psi$ maps $C(a(\mathcal{T}))$ onto $C(a(\mathcal{W})) \oplus \IC.$  If $(f,c) \in C(a(\mathcal{W})) \oplus \IC,$ then define $g: a(\mathcal{T}) \rightarrow \IC$ by \begin{equation*} g(\{\tilde{x}_F\}) = \left\{ \begin{array}{ll} f(\{\tilde{x}_F\}), & \{\tilde{x}_F\}_F \neq \{(0, \ldots 0)_F\} \\ c, & \{\tilde{x}_F\}=\{(0, \ldots, 0)_F\}\end{array} \right. . \end{equation*}
    If $V$ is an open set in $\IC$, then $g^{-1}(V)$ is equal to either $f^{-1}(V)$ or $f^{-1}(V) \cup \{(0, \ldots,0)_F\}$, both of which are open sets in $a(\mathcal{T})$ since $f \in C(a(\mathcal{W})).$  
    So $g \in C(a(\mathcal{T}))$, and $\psi$ is surjective.
 \end{proof}
 
 We now combine our results to prove our main theorem:
 
 \begin{theorem}  
For $0 < s <1,$ let $\varphi_{s}(z)=sz+(1-s).$  Let $\mathcal{C}$ denote the commutator ideal of $C^*\left(\left\{C_{\varphi_s}: 0 < s < 1\right\}\right).$  Then there exists a $*$-homomorphism $\psi :C^*\left(\left\{C_{\varphi_s}: 0 < s < 1\right\}\right) \rightarrow  AP(\IR) \oplus \IC$ such that    $$ 0 \rightarrow \mathcal{C}  \hookrightarrow C^*\left(\left\{C_{\varphi_s} :
   0 < s < 1\right\} \right) \stackrel{\psi}{\rightarrow} AP(\mathbb{R}) \oplus \mathbb{C} \rightarrow 0$$ is a short exact sequence.
 \end{theorem}
 
 \begin{proof} Let $a(\mathcal{W})$ and $a(\mathcal{T})$ be defined as above.  Since $a(\mathcal{W})$ and  $a(\mathcal{T})$ are inverse limits of $\{\{\sigma_F\}, \{\pi_F\}\}$ and $\{\{\tilde{\sigma}_F\}, \{\tilde{\pi}_F\}\},$ respectively, $a(\mathcal{W})$ and $a(\mathcal{T})$ are homeomorphic to $\sigma_{ap}(\left\{W_{\chi_{\alpha}}: \alpha \in \IR^+\right\})$ and $\sigma_{ap}(\left\{e^{-\frac{\alpha}{2}}T_{e^{-\alpha z}} : \alpha \in \IR^+\right\})$, respectively.
   Thus, by Corollary \ref{contiso}, \begin{equation*}\label{isoweneed} C\left(\sigma_{ap}\left(\left\{e^{-\frac{\alpha}{2}}T_{e^{-\alpha z}} : \alpha \in \IR^+\right\}\right)\right) \cong C\left(\sigma_{ap}\left(\left\{W_{\chi_{\alpha}}: \alpha \in \IR^+\right\}\right)\right) \oplus \IC.\end{equation*}
 
 Let $\mathcal{C_W}$ and $\mathcal{C_T}$ denote the commutator ideals of $C^*\left(\left\{W_{\chi_{\alpha}}: \alpha \in \IR^+\right\}\right)$ and   $C^*\left(\left\{e^{-\frac{\alpha}{2}}T_{e^{-\alpha z}} : \alpha \in \IR^+\right\}\right),$ respectively.  Then by Theorem \ref{Bunceinfinite}, \begin{align*}\label{intequiv} C^*\left(\left\{e^{-\frac{\alpha}{2}}T_{e^{-\alpha z}} : \alpha \in \IR^+\right\}\right) / \mathcal{C_T} &\cong  C\left(\sigma_{ap}\left(\left\{e^{-\frac{\alpha}{2}}T_{e^{-\alpha z}} : \alpha \in \IR^+\right\}\right)\right) \\ &\cong \left(C^*\left(\left\{W_{\chi_{\alpha}}: \alpha \in \IR^+\right\}\right) / \mathcal{C_W}\right) \oplus \IC.\end{align*}
 Applying Theorem \ref{CDiso}, we obtain
 \begin{equation*}\label{almostthere} 
 C^*(\left\{e^{-\frac{\alpha}{2}}T_{e^{-\alpha z}} : \alpha \in \IR^+\right\}) / \mathcal{C_T} \cong AP(\IR) \oplus \IC. \end{equation*}  
 By Theorem \ref{cowenkrieteunitary}, $$C^*\left(\left\{C_{\varphi_s}: 0 < s < 1\right\}\right)/ \mathcal{C} \cong C^*(\left\{e^{-\frac{\alpha}{2}}T_{e^{-\alpha z}} : \alpha \in \IR^+\right\}) / \mathcal{C_T},$$ which yields the desired result.
 \end{proof}

We can explicitly describe how the $*$-homomorphism $\psi$ acts on the generators of $C^*\left(\left\{C_{\varphi_s}: 0 < s < 1\right\}\right)$.  For $0 < s <1$, $\psi(C_{\varphi_s})=\left(s^{-\frac{1}{2}+iy}, 0\right),$ and $\psi(I)=(1,1).$  Using this description, we can obtain spectral information for a dense set of operators in $C^*\left(\left\{C_{\varphi_s}: 0 < s < 1\right\}\right)$.  To simplify notation in the following result, we extend the definition of $C_{\varphi_s}$ to include $C_{\varphi_1}=I.$  Then every word in the generators of $C^*(\{C_{\varphi_s} : 0 < s< 1\})$ can be written in the form \begin{equation}\label{form} C_{\varphi_{s_1}} C_{\varphi_{s_2}}^* C_{\varphi_{s_3}} \ldots C_{\varphi_{s_m}}^*, \end{equation} where $m$ is a positive, even integer, $s_1, s_m \in (0,1],$ and $s_2, s_3, \ldots s_{m-1} \in (0,1).$  Note that a word of form (\ref{form}) is the identity operator if and only if $m=2$ and $s_1=s_2=1.$

\begin{corollary} Let $n \in \IN$ and $c_0, c_1, \ldots, c_n \in \IC.$  For each $i \in \{1, \ldots, n\},$ let $m_i$ be a positive, even integer, and let $s_{i,2}, s_{i,3}, \ldots, s_{i,m_{i-1}} \in (0,1)$ and $s_{i,1}, s_{i, m_i} \in (0,1],$ with either $s_{i,1} \neq 1$ or $s_{i,m_i} \neq 1$ if $m_i=2.$  Consider the operator \begin{equation*} A=c_0I + \sum_{i=1}^nc_i C_{\varphi_{s_{i,1}}}C_{\varphi_{s_{i, 2}}}^*C_{\varphi_{s_{i,3}}} \ldots C_{\varphi_{s_{i,m_i}}}^* \in C^*\left(\left\{C_{\varphi_s}: 0 < s < 1\right\}\right).\end{equation*} Then 
\begin{equation*} \{c_0\} \cup \overline{\left\{c_0 +\sum_{i=1}^{n} c_i s_{i,1}^{-\frac{1}{2}+iy} s_{i,2}^{-\frac{1}{2}-iy} s_{i,3}^{-\frac{1}{2}+iy} \ldots s_{i, m_i}^{-\frac{1}{2}-iy} : y \in \IR\right\}} \subseteq \sigma(A).\end{equation*} \end{corollary}

\section{\textsc{Is} $C^*\left(\left\{C_{\varphi_s}: 0 < s < 1\right\}\right)$ \textsc{Irreducible?}}

In this section, we consider the C$^*$-algebra $C^*\left(\left\{\tilde{T}_{s^{\frac{z}{1-z}}} : 0 < s < 1\right\}\right)$ of operators on $H^2(\tilde{\mu})$, which is unitarily equivalent to $C^*\left(\left\{C_{\varphi_s}: 0 < s < 1\right\}\right)$.  The main reason for using this setting is that $\tilde{\mu}$ is a compactly supported measure on $\IC$, which allows us to take advantage of some known facts about subnormal operators.  We begin by recalling these facts.

Let $\nu$ be a compactly supported measure on $\IC.$  Let $P^2(\nu)$  denote the closure of the polynomials in $L^2(\nu).$  We define the operator $S_{\nu} : P^2(\nu) \rightarrow P^2(\nu)$ by $(S_{\nu} f)(z) = zf(z)$ for all $f \in P^2(\nu).$  Note that $S_{\nu}$ is a subnormal operator.  Let $\{S_{\nu}\}^{\prime}$ denote the commutant of $S_{\nu}.$ Then a corollary to a theorem of T. Yoshino \cite{Yoshino:1969} states:
\begin{theorem}\label{Yoshino} Let $\nu$ be a compactly supported measure on $\IC$. Then $$\{S_{\nu}\}^{\prime} = \left\{M_{\psi}|_{P^2(\mu)} : \psi \in P^2(\nu) \cap L^{\infty}(\nu)\right\},$$ where $M_{\psi}$ denotes the multiplication operator $M_{\psi}f=\psi f$ on $L^2(\IC, \nu).$ \end{theorem} This corollary and other facts about $P^2(\nu)$ and $S_{\nu}$ can be found in \cite{Conway:1991}.
For our investigations, we will take $\nu = \tilde{\mu}$.  Note that $P^2(\tilde{\mu})=H^2(\tilde{\mu})$, and   for any $\psi \in P^2(\tilde{\mu}) \cap L^{\infty}(\tilde{\mu})$, $M_{\psi}|_{P^2(\tilde{\mu})} = \tilde{T}_{\psi}$.  Moreover, $S_{\tilde{\mu}} = \tilde{T}_z.$

Before we address the question of irreducibility, we establish a lemma about functions in $H^2(\tilde{\mu})$ that will play a key role in the proof.

\begin{lemma} \label{constant} If $f \in H^2(\tilde{\mu})$, then $P_{\tilde{\mu}}\overline{f}$ is a constant function. \end{lemma}

\begin{proof}

First, let $0 < a, b \leq 1.$  Since $U_1U_2^*$ is a unitary operator from $H^2(\tilde{\mu})$ onto $\mathcal{N}$ with $U_1U_2^*(ab)^{\frac{z}{1-z}}= (ab)^z$, we have that \begin{align} \left<a^{\frac{z}{1-z}}, P_{\tilde{\mu}}\overline{b^{\frac{z}{1-z}}} \right>_{H^2(\tilde{\mu})} &= \left<a^{\frac{z}{1-z}}, \overline{b^{\frac{z}{1-z}}}\right>_{L^2(\tilde{\mu})} \notag \\
&= \left<(ab)^{\frac{z}{1-z}}, 1 \right>_{H^2(\tilde{\mu})} \notag \\
&= \left<(ab)^z, 1 \right>_{\mathcal{N}} \notag \\
&= (ab)^0=1 \label{innerprod} \end{align} 
since $1=K_0 \in \mathcal{N}.$  Thus, $$\left<a^{\frac{z}{1-z}}, P_{\tilde{\mu}}\overline{b^{\frac{z}{1-z}}} -1 \right>_{H^2(\tilde{\mu})} = \left<a^{\frac{z}{1-z}}, P_{\tilde{\mu}}\overline{b^{\frac{z}{1-z}}} -P_{\tilde{\mu}}\overline{1^{\frac{z}{1-z}}} \right>_{H^2(\tilde{\mu})} = 0.$$  
Since $0 < a \leq 1$ was arbitrary and the linear span of $\left\{a^{\frac{z}{1-z}} : 0 < a \leq 1\right\}$ is dense in $H^2(\tilde{\mu})$, $P_{\tilde{\mu}} \overline{b^{\frac{z}{1-z}}} = 1.$

Let $f \in H^2(\tilde{\mu}).$  Then there is a sequence of the form $\left\{\sum_{j=1}^{N_n} c_{n,j} b_{n,j}^{\frac{z}{1-z}} \right\}_{n=1}^{\infty}$ that converges to $f$ in $L^2(\tilde{\mu}).$  
 Here, each $c_{n,j} \in \IC,$ and each $b_{n,j} \in (0,1]$.
 
 Set $\gamma := \left<f, 1 \right>_{H^2(\tilde{\mu})},$
and  
let $\varepsilon >0$ be given. Then there exists $M>0$ such that if $n >M$, then $\left|\left|\sum_{j=1}^{N_n} c_{n,j} b_{n,j}^{\frac{z}{1-z}} -f\right|\right|_{L^2(\tilde{\mu})} < \frac{\varepsilon}{2}.$  Thus, for $n>M$,
 \begin{equation*} \left|\overline{\gamma} - \sum_{j=1}^{N_n}\overline{c_{n,j}}\right| = \left| \overline{\left<f- \sum_{j=1}^{N_n} c_{n,j} b_{n,j}^{\frac{z}{1-z}},1\right>_{H^2(\tilde{\mu})}} \right|  \leq \left| \left| f - \sum_{j=1}^{N_n} c_{n,j} b_{n,j}^{\frac{z}{1-z}} \right| \right|_{L^2(\tilde{\mu})} <  \frac{\varepsilon}{2} \end{equation*} by (\ref{innerprod}) and the fact that $\tilde{\mu}(\overline{\ID})=1.$  Also, for $n >M$, \begin{align*} \left| \left|P_{\tilde{\mu}} \overline{f} - \sum_{j=1}^{N_n} \overline{c_{n,j}} \right| \right|_{H^2(\tilde{\mu})} &= \left| \left|P_{\tilde{\mu}} \overline{f} - P_{\tilde{\mu}}\overline{\sum_{j=1}^{N_n} c_{n,j} b_{n,j}^{\frac{z}{1-z}} } \right| \right|_{H^2(\tilde{\mu})} \\
&\leq \left|\left|f - \sum_{j=1}^{N_n} c_{n,j} b_{n,j}^{\frac{z}{1-z}}  \right| \right|_{L^2(\tilde{\mu})}< \frac{\varepsilon}{2}.\end{align*} Hence $\left|\left| P_{\tilde{\mu}}\overline{f}- \overline{\gamma}\right| \right|_{H^2(\tilde{\mu})} < \varepsilon.$  Since $\varepsilon$ was arbitrary, $\left(P_{\tilde{\mu}}\overline{f}\right)(z)=\overline{\gamma}$ $\tilde{\mu}$-almost everywhere.
\end{proof}

We now return to the question of whether $C^*\left(\left\{C_{\varphi_s} : 0 < s < 1\right\}\right)$ is irreducible.  We will show that the C$^*$-algebra is irreducible by showing that the commutant of $C^*\left(\left\{\tilde{T}_{s^{\frac{z}{1-z}}} : 0 < s <1\right\}\right)$ consists of only scalar multiples of the identity operator on $H^2(\tilde{\mu})$.

\begin{theorem} The C$^*$-algebra $C^*\left(\left\{C_{\varphi_s} : 0 < s < 1\right\}\right)$ is irreducible.
\end{theorem}
\begin{proof}
       
For any $f, g \in H^2(\tilde{\mu})$, we calculate
\begin{align} \int_{0}^1 \left< \tilde{T}_{s^{\frac{z}{1-z}}} f, g \right>_{H^2(\tilde{\mu})} ds &= \int_{0}^1 \int_{\overline{\ID}} s^{\frac{z}{1-z}} f(z) \overline{g(z)} d\tilde{\mu}(z) ds \notag \\
&= \int_{0}^{\infty} \int_{\overline{\ID}} e^{-t} e^{-t \frac{z}{1-z}} f(z) \overline{g(z)} d\tilde{\mu}(z) dt \label{change} \\
&= \int_{\overline{\ID}} \int_0^{\infty} e^{\frac{-t}{1-z}} f(z)\overline{g(z)} dt d\tilde{\mu}(z) \label{Fubini}\\
&= \int_{\overline{\ID}} (1-z)f(z)\overline{g(z)} d\tilde{\mu}(z) \label{nomass}\\
&= \left<\tilde{T}_{1-z} f, g \right>_{H^2(\tilde{\mu})}, \notag
\end{align}
where (\ref{change}) is obtained via the change of variables $s=e^{-t}$.  The use of Fubini's theorem in (\ref{Fubini}) can be justified by the fact that $\left|e^{-\frac{t}{1-z}}\right| \leq e^{-\frac{t}{2}}$ for all $t \geq 0$ and $z \in \overline{\ID} \setminus \{1\}$.  The results in  (\ref{Fubini}) and $(\ref{nomass})$ both rely on the fact that $\mu(\{1\})=0.$

Suppose $A$ is in the commutant of $C^{*}\left(\left\{\tilde{T}_{s^{\frac{z}{1-z}}} : 0 < s < 1\right\}\right)$.  
Since $A$ commutes with  $\tilde{T}_{s^{\frac{z}{1-z}}}$ for all $0 < s <1$, then, for all $f, g \in H^2(\tilde{\mu})$, 
\begin{align*} \left<A \tilde{T}_{1-z} f, g \right>_{H^2(\tilde{\mu})} & = \left<\tilde{T}_{1-z}f, A^{*}g\right>_{H^2(\tilde{\mu})}\\
&= \int_0^1 \left<\tilde{T}_{s^{\frac{z}{1-z}}}f, A^*g\right>_{H^2(\tilde{\mu})} ds\\
& = \int_0^1\left<\tilde{T}_{s^{\frac{z}{1-z}}}Af, g \right>_{H^2(\tilde{\mu})} ds\\
& = \left<\tilde{T}_{1-z}Af, g\right>_{H^2(\tilde{\mu})}.
\end{align*}
Thus, $A$ commutes with $\tilde{T}_{1-z}$ and hence $\tilde{T}_{z}$.  Since $A$ must also commute with $\tilde{T}^*_{s^{\frac{z}{1-z}}}$ for all $0 < s <1$, $A^*$ also commutes with $\tilde{T}_z.$

Then, by Theorem \ref{Yoshino}, there exists $\psi, \rho \in H^2(\tilde{\mu}) \cap L^{\infty}(\tilde{\mu})$ such that $A=\tilde{T}_{\psi}$ and $A^*=\tilde{T}_{\rho}.$  Since $\tilde{T}_{\rho}^*=\tilde{T}_{\overline{\rho}}$, we have that $\tilde{T}_{\psi}=\tilde{T}_{\overline{\rho}}$ and $$||\psi-P_{\tilde{\mu}} \overline{\rho}||_{H^2(\tilde{\mu})}=||(\tilde{T}_{\psi}-\tilde{T}_{\overline{\rho}})1||_{H^2(\tilde{\mu})} =0.$$  But, by Lemma \ref{constant}, there exists $\gamma \in \IC$ such that $\left(P_{\tilde{\mu}}\overline{\rho}\right)(z) = \gamma$ $\tilde{\mu}$-almost everywhere.  Therefore, $\psi(z) = \gamma$ $\tilde{\mu}$-almost everywhere, and $A=\tilde{T}_{\psi}=\gamma I.$ \end{proof}
\vspace{12pt}

\textit{Acknowledgement.} The author thanks her advisor, Thomas Kriete, for his guidance, encouragement, and many helpful conversations throughout the completion of this work.

\end{document}